\documentclass{amsart}

\usepackage[latin1]{inputenc}

\usepackage{amsmath,amssymb,amscd,latexsym,color}
\usepackage[dvips]{graphicx}

\usepackage[english]{babel}

\DeclareMathOperator{\Cyl}{Cyl}
\DeclareMathOperator{\Shadow}{Shadow}
\DeclareMathOperator{\Shade}{Shade}

\DeclareMathOperator{\supp}{supp}

\DeclareMathOperator{\Coex}{Coex}

\def\shell{\text{Shell}}

\def\gstoch{\succ}

\newcommand{\chapeau}[1]{\overrightarrow{#1}}

\newcommand{\Peierls}{\mathcal{P}\text{eierls}}

\newcommand{\N}{\mathbb N}

\newcommand{\Z}{\mathbb{Z}}
\newcommand{\Zdeux}{\mathbb{Z}^2}
\newcommand{\Zdeuxstar}{\mathbb{Z}_{*}^2} 
\newcommand{\Zp}{\mathbb{Z}_{*}^2}
\newcommand{\Zd}{\mathbb{Z}^d}

\newcommand{\R}{\mathbb{R}}
\newcommand{\Rd}{\mathbb{R}^d}

\renewcommand{\P}{\mathbb{P}}
\newcommand{\E}{\mathbb{E}\ }

\newcommand{\Ld}{\mathbb{L}^d}
\newcommand{\Ldeux}{\mathbb{L}^2}
\newcommand{\Ldeuxstar}{\mathbb{L}_{*}^2}

\newcommand{\Ed}{\mathbb{E}^d}
\newcommand{\Edeux}{\mathbb{E}^2}
\newcommand{\Edeuxstar}{\mathbb{E}^2_{*}}

\renewcommand{\epsilon}{\varepsilon}
\renewcommand{\phi}{\varphi}
\renewcommand{\limsup}{\overline{\lim}}
\renewcommand{\liminf}{\underline{\lim}}

\newcommand{\ie}{\emph{i.e. }}

\newcommand{\resp}{\emph{resp. }}
\newcommand{\vsp}{\vspace{0.2cm}}

\newcommand{\miniop}[3]{%
\renewcommand{\arraystretch}{0.6}
\begin{array}{c}
{\scriptstyle #1}\\
#2\\
{\scriptstyle #3}
\end{array}
\renewcommand{\arraystretch}{1}}

\newcommand{\Card}[1]{\vert #1 \vert}

\newcommand{\norme}[2]{\|#1 \|_{#2}}

\newcommand{\1}{1\hspace{-1.3mm}1}

\newcommand{\pcfleche}{\overrightarrow{p_c}}

\newcommand{\froufrou}{\partial_{\text{ext}}}

\begin{document}

{
\newtheorem{theorem}{Theorem}[section]
\newtheorem{conjecture}[theorem]{Conjecture}

}
\newtheorem{lemme}[theorem]{Lemma}
\newtheorem{defi}[theorem]{Definition}
\newtheorem{coro}[theorem]{Corollary}
\newtheorem{rem}[theorem]{Remark}
\newtheorem{prop}[theorem]{Proposition}

\title[First-passage competition and density]{First-passage competition with different speeds: positive density for both species  is impossible}

{
\author{Olivier Garet}
\author{R{\'e}gine Marchand}
\address{Laboratoire de Math{\'e}matiques, Applications et Physique
Math{\'e}matique d'Orl{\'e}ans UMR 6628\\ Universit{\'e} d'Orl{\'e}ans\\ B.P.
6759\\
 45067 Orl{\'e}ans Cedex 2 France}
\email{Olivier.Garet@univ-orleans.fr}

\address{Institut Elie Cartan Nancy (math{\'e}matiques)\\
Universit{\'e} Henri Poincar{\'e} Nancy 1\\
Campus Scientifique, BP 239 \\
54506 Vandoeuvre-l{\`e}s-Nancy  Cedex France}
\email{Regine.Marchand@iecn.u-nancy.fr}
}

\def\motsclefs{First-passage percolation, competition, coexistence, random growth, moderate deviations.}


\begin{abstract}
Consider two epidemics whose expansions on $\Zd$ are governed by two families of passage times that are distinct and stochastically comparable. We prove that when the weak infection survives, the space occupied by the strong one is almost impossible to detect: for instance, it could not be observed by a medium resolution satellite. We also recover the same fluctuations with respect to the asymptotic shape as in the case where the weak infection evolves alone. In dimension two, we prove that one species finally occupies a set with full density, while the other one only occupies a set of null density. We also prove that the Häggström-Pemantle non-coexistence result
"except perhaps for a denumerable set" can be extended to families
of stochastically comparable passage times indexed by a continuous parameter.
\end{abstract}

{\maketitle 
}
\setcounter{tocdepth}{1}

\section{Introduction}

Consider two species both trying to colonize the graph $\Zd$. The expansion of each species is governed by independent identically distributed random passage times
attached to the bonds of the graph, as in first-passage percolation, and each vertex of the graph can only be infected once, by the first species that reaches it. Is it possible that both species simultaneously succeed in invading an infinite
subset of the net, in other words that coexistence occurs?
That is the kind of question which was asked in the middle of the 90's
by Häggström and Pemantle in two seminal papers~\cite{Haggstrom-Pemantle-1,
Haggstrom-Pemantle-2}, where they gave the
first results towards the following conjectures:
\begin{itemize}
\item If the two species travel at the same speed, coexistence is possible.
\item If one of them travels faster than the other one, coexistence is impossible.
\end{itemize}
The passage times considered by 
Häggström and Pemantle follow exponential laws, which provides a Markov property and allows a description of the competition process in terms of particle system. However, the first-passage percolation setting naturally enables to consider competition with general passage times, even if the Markovianity is lost.

The problem of coexistence for two similar species has been solved
by Häggström and Pemantle~\cite{Haggstrom-Pemantle-1} in dimension two for exponential passage times, then by the authors of the present paper in any dimension for general passage times, under assumptions that are close
to optimality~\cite{GM-coex}. Shortly later, Hoffman~\cite{hoffman} gave a different proof involving
tools that seem to allow an extension to a larger number of species --
see Hoffman's manuscript~\cite{hoffman-geo}.

On the contrary, the state of the art about the second conjecture -- the non-coexistence problem -- did not much change since its statement. More precisely, if one species travels according passage times
following the exponential law with intensity $1$, while the other one 
 travels according passage times
following the exponential law with  intensity $\lambda\ne 1$, it is believed
that coexistence is not possible. However, it is only known that coexistence
is not possible "except perhaps for a denumerable set of values of $\lambda$", as it was proved by Häggström and Pemantle~\cite{Haggstrom-Pemantle-2}. 
To sum up, if one denotes by $\text{Coex}$ the set of intensities for the second particle that allow coexistence, we know that $\text{Coex}\supset\{1\}$
and $\text{Coex}$ is denumerable, but we would like to have $\text{Coex}=\{1\}$.
It follows that we are currently in the following perplexing situation: we
know that for almost every value of $\lambda$, coexistence does not happen, 
but we are unable to exhibit any value of $\lambda$ such that coexistence does 
not occur.

Therefore, the aim of this paper is to prove a weakened version of non-coexistence for epidemics with distinct speeds. Let us first introduce our framework: we consider two epidemics whose expansions are governed by two families of independent and identically distributed passage times whose laws are distinct and stochastically comparable, which of course includes the case of exponential laws. We say that \textit{strong coexistence} occurs when each species finally occupy a set with positive natural density. 

In dimension two, we prove that, almost surely, strong coexistence does not occur. More precisely, we show that almost surely, at infinite time, one species fills a set with full natural density, whence the
other one only  fills a set with null natural density. In higher dimension, connectivity problems prevent us to obtain such a complete result. However, we show that, roughly speaking, a medium resolution satellite only sees one type of particles. 

By the way, we also prove that the Häggström-Pemantle non-coexistence result
"except perhaps for a denumerable set" can be extended to families
of stochastically comparable passage times indexed by a continuous parameter.
Note that the Häggström-Pemantle method~\cite{Haggstrom-Pemantle-2} to prove denumerability of 
 $\text{Coex}$ has already been transposed to other
models having familiarities with first-passage percolation: at first
by Deijfen, Häggström, and Bagley for a model with spherical symmetry~\cite{DHB}, then by the
authors of the present paper for some percolating model~\cite{GM-bernoulli}.

Before giving more rigorous statements of our results, let us introduce general notations and give a formal description of the competition model.

\subsection*{General notations}
We denote by $\Z$ the set of integers, by $\N$ the set of non negative integers.

We endow the set $\Zd$ with the set of edges $\Ed$ between sites of $\Zd$ that are at distance $1$ for the Euclidean distance: the obtained graph is denoted by $\Ld$. Two sites $x$ and $y$ that are linked by an edge are said to be \emph{neighbors} and this relation is denoted: $x \sim y$. If $A$ is a subset of $\Zd$, we define the border of $A$:
$$\partial A=\{ z \in A: \quad \exists y \in A^c \quad y \sim z\}.$$
A path in $\Zd$ is a sequence $x_0,x_1, \dots, x_l$ of points in $\Zd$ such that two successive points are neighbors. The integer $l$ is called the length of the path.

The critical percolation for Bernoulli percolation (oriented percolation) on $\Zd$ is denoted by $p_c=p_c(d)$ (respectively,  $\pcfleche=\pcfleche(d)$).

Let us now recall the concept of stochastic domination: we say that a probability 
measure $\mu$ dominates a probability measure $\nu$, which is denoted by $\nu\prec\mu$, if 
$${\int f\ d\nu}\le{\int f\ d\mu}$$ 
holds as soon as $f$ in a non decreasing function. 

The complementary event of $A$ will mostly be denoted $A^c$. But sometimes, to improve readability, we prefer to use $\complement A$.

\subsection*{Assumptions on passage times}
Let $\nu_{p_1}$ and $\nu_{p_2}$
be two probability measures on $[0,+\infty)$. We will always assume that
\begin{eqnarray*}
(H1) &&  \nu_{p_1}\gstoch \nu_{p_2} \text{ and } \nu_{p_1} \neq \nu_{p_2}. \\
(H2) && \forall k,l \in \mathbb{N} \quad \nu_{p_1}^{*k}\otimes \nu_{p_2}^{*l}(\{(x,x): \, x \in [0,+\infty)\})=0.\\
(H3) && \text{for } i \in \{1,2\}, \; \nu_{p_i}(0)<p_c. \\
(H4) && \text{for } i \in \{1,2\}, \; \nu_{p_i}(\inf \supp \nu_{p_i}) < \pcfleche(d).\\
(H5) && \text{for } i \in \{1,2\}, \; \exists \gamma>0 \text{ such that }  \int_0^{+\infty} \exp(\gamma x)\ d\nu_{p_i}(x)< +\infty.
\end{eqnarray*}
In Assumption $(H_3)$, $\inf \supp \nu$ denotes the infimum of the support of the measure $\nu$. Note that Assumptions $(H2)$, $(H3)$, and $(H4)$ are clearly fulfilled when $\nu_{p_1}$ and $\nu_{p_2}$ are continuous with respect to Lebesgue's measure.

Assumption $(H2)$ exactly says that a sum of $k$ independent random variables with common distribution $\nu_{p_1}$ and a sum of $l$ independent random variables with common distribution $\nu_{p_2}$, independent of the first family, have probability $0$ to be equal: this will ensure that, during the competition process, no vertex of $\Zd$ can be reached exactly at the same time by the two epidemics. 
Assumptions $(H_1)$, $(H3)$ and $(H_4)$ are the ones used by van den Berg and Kesten in~\cite{vdb-kes} to prove that, in first-passage percolation, the time constant for $\nu_{p_2}$ is strictly smaller than the one for $\nu_{p_1}$.
Finally, assumption $(H5)$ gives access to large deviations and moderate deviations related to the asymptotic shape in first-passage percolation.

\subsection*{Construction of the competition model}
The first infection (second infection) will use independent identically distributed  passage times with common law $\nu_{p_1}$ (respectively,  $\nu_{p_2}$) and start from the source $s_1$ (respectively,  $s_2$, distinct from $s_1$).
As $\nu_{p_1} \gstoch \nu_{p_2}$, species $1$ will be slower (or weaker)  than species $2$. 

First, we couple the two measures $\nu_{p_1}$ and $\nu_{p_2}$ in agreement with the stochastic comparison relation $(H_1)$:
there exists a probability measure  $m$  on $[0,+\infty)\times[0,+\infty)$ such that
$m(\{(x,y)\in [0,+\infty); x\ge y\})=1$ and the marginals of $m$ are $\nu_{p_1}$ and $\nu_{p_2}$.

Now, we consider, on $\Omega=([0,+\infty)\times[0,+\infty))^{\Ed}$, the probability measure $\P=m^{\otimes\Ed}$. For $\omega=(\omega^1_e, \omega^2_e)_{e \in \Ed}$, the number $\omega^i_e$ represents the time needed by species $i$ to cross edge $e$. Note that under $\P$, 
for each $i\in\{1,2\}$, the variables $(\omega^i_e)_{e\in\Ed}$  are independent identically distributed  with common law $\nu_{p_i}$.
Moreover, we almost surely have
$$\forall e\in\Ed \quad \omega^1_e \ge \omega^2_e.$$
It remains to construct the two infections in a realization $\omega \in \Omega$. 
Let  $E=([0,+\infty]\times [0,+\infty])^{\Zd}$. 
We recursively define a $E$-valued sequence $(X_n)_{n\ge 0}$ and a non-negative
sequence  $(T_n)_{n\ge 0}$. The sequence $(T_n)_{n\ge 0}$ contains the successive times of infections, while a point $\epsilon=(\epsilon^1(z), \epsilon^2(z))_{z \in \Zd} \in E$ codes, for each site $z$, its times of infection $\epsilon^1(z)$ (its time of infection $\epsilon^2(z)$) by the first infection (respectively,  second infection). We start the process with two distinct sources $s_1$ and $s_2$ in $\Zd$, and set $T_0=0$ and 
$$
X_0=(X_0^1(z), X_0^2(z))_{z \in \Zd} \text{ with } X_0^i(z)=0 \text{ if } z=s_i, \text{ and } X_0^i(z)=+\infty \text{ otherwise.}
$$
This means that at time $0$, no point of $\Zd$ has yet been infected but the two initial sources $s_1$ and $s_2$. Then, for $n \ge 0$, define the next time of infection:
$$T_{n+1}=\inf\{ X_n^i(y)+\omega_{\{y,z\}}^i: \; \{y,z\}\in\Ed, \; i\in\{1,2\}, \;  X_n^{3-i}(z)=+\infty\}.$$
Note that 
Assumption $(H2)$ ensures that if this infimum is reached for several triplets $(i_1,y_1,z_1)$, \dots, $(i_l,y_l,z_l)$, then $i_1=\dots=i_l=i$. In this case, the next infections are of type $i$ from each $y_j$ to each $z_j$. The set of infected points of type $3-i$ has not changed:
$$\forall x \in \Zd \quad X_{n+1}^{3-i}(x)=X_n^{3-i}(x),$$
while the points $z_j$ has been infected by species $i$ at time $X_n^i(y_j)+\omega_{\{y_j,z_j\}}^i$:
$$
\begin{array}{ll}
& \forall x \in \Zd\backslash \{z_1, \dots, z_l\} \quad  X_{n+1}^{i}(x)=X_n^{i}(x)\\
 \text{ and } & \forall j \in \{1, \dots, l\} \quad X_{n+1}^i(z_j)=X_n^i(y_j)+\omega_{\{y_j,z_j\}}^i.
\end{array}$$
Note that $X_n^i(y)$ and $X_n^{3-i}(y)$ can not be simultaneously finite, which corresponds to the fact that each site is infected by at most one type of infection. Moreover, once $\min(X_n^1(x), X_n^2(x))$ is finite, its value -- the time of infection of $x$ -- does not change any more with $n$.

Note however that this coupling is nothing else than a useful tool for our proofs: this does note constrain the evolution of the process. Particularly, the very definition of the evolution process tells us that the $(\omega^e_i)_{e \in \Ed, i \in \{1,2\}}$ could be independent as well without the law of the evolution process being changed.

We can now define the sets $\eta(t),\eta^1(t),\eta^2(t)$ that are respectively the sets of infected points, infected points of type $1$, infected points of type $2$ at time $t$, by setting
$$\forall i\in\{1,2\}\quad \forall t\in [T_n,T_{n+1})\quad  \eta^i(t)=\{z\in\Zd: \; X_n^i(z)<+\infty\}$$
and $\eta(t)=\eta^1(t)\cup \eta^2(t)$. 
We also introduce the the sets of points that are finally infected by each epidemic:
$$\forall i\in\{1,2\}\quad  \eta^i(\infty)=\miniop{}{\cup}{t\ge 0}\eta^i(t).$$
The set
$
\mathcal{G}^i = \left\{ |\eta^i(\infty)|=+\infty\right\} \text{ for } i=1,2,
$
corresponds to the unbounded growth of species $i$, and coexistence is thus the event $\Coex=\mathcal{G}^1 \cap \mathcal{G}^2.$

\subsection*{Coupling with first-passage percolation}
The evolution of the two infections can be compared with the evolution of classical first-passage percolation.

Assume that $\nu$ is a probability measure on $[0, +\infty)$, such that
$$\nu(0)<p_c \text{ and } 
\exists \gamma>0 \text{ such that } \int_0^{+\infty}\exp(\gamma x)\ d\nu(x)< +\infty.
$$
On $\Omega=[0,+\infty)^{\Ed}$, consider the probability measure $\P=\nu^{\otimes\Ed}$, which makes the coordinates  $(\omega_e)_{e\in\Ed}$  independent identically distributed  random variables with common law $\nu$. Then for each $x\in\Zd$ and $t\ge 0$, we define the set of points reached from $x$ in a time less than $t$:
$$B^{x}(t)=\left\{y\in\Zd: \text{ there exists a path }\gamma \text{ from }x\text{ to }y\text{, with }\sum_{e\in \gamma} \omega_e\le t\right\}.$$
The classical shape theorem gives the existence of 
a norm $\|.\|_{\nu}$ on $\Rd$ such that 
$B^0(t)/t$ almost surely converges to the unit ball $\mathcal B$ for $\|.\|_{\nu}$. 

Note that the competition model contains two simple first-passage percolation models: 
for each $i\in\{1,2\}$, $x\in\Zd$ and $t\ge 0$, the set
$$B^{x}_{p_i}(t)=\left\{y\in\Zd: \text{ there exists a path }\gamma \text{ from }x\text{ to }y\text{, with }\sum_{e\in \gamma} \omega_e^i\le t\right\}$$
is the random ball of radius $t$ of first-passage 
percolation starting from $x$ with passage time law $\nu_{p_i}$. 
For simplicity, the related norm will be denoted by $\norme{.}{p_i}$, and its associated discrete balls $\mathcal B_{p_i}^x(t)=\{y \in \Zd: \; \norme{y-x}{p_i} \le t\}$.
 
Here are the coupling relations between the competition model and first-passage percolation:
\begin{lemme}
\label{coupling}
$\forall t>0 \quad  \eta^1(t) \subset B_{p_1}^{s_1}(t) \text{ and } 
\eta^2(t)  \subset  B_{p_2}^{s_2}(t)   \text{ and } B_{p_1}^{s_1}(t) \subset \eta(t).$
\end{lemme}
We postpone the proof of these (not so) obvious properties to the next section.

\subsection*{Statement of results} 
We denote $\|.\|_2$ the euclidean norm on $\Rd$, by $\langle .,.\rangle$ the corresponding scalar product and by  $\mathcal S$ the corresponding unit sphere: $\mathcal S  =  \{x\in\R^2 : \; \|x\|_{2}=1\}$.
Let $y,z\in \Rd$, $\chapeau x\in \mathcal S$, and $R,h>0$. We define:
\begin{align*}
d(y,\R\chapeau x) & =  \| y-\langle y,\chapeau x\rangle\chapeau x\|_2\text{ (the Euclidean distance from $y$ to the line $\R\chapeau x$)}, \\
\Cyl_z(\chapeau x,R,h) & =  \{y\in\Zd : \; d(y-z,\R\chapeau x)\le R\text{ and }0\le \langle y-z,\chapeau x \rangle \le h\}, \\
\Cyl_+(\chapeau x ,R) & =  \Cyl_0(\chapeau x,R,+\infty),\\
\Cyl(\chapeau x,R) & =  \{y\in\Zd : \;  d_2(y,\R\chapeau x )\le R\}.
\end{align*}
We can then define the following events:
\begin{align*}
\Shadow(\chapeau x,t, R) &  =  \left\{\partial\eta(t)\cap\eta^2(t)\text{ disconnects }\eta^1(t)\text{ from infinity in }\Cyl_+(\chapeau x ,R)\right\}, \\
\Shade(t, R) & =  \miniop{}{\cup}{x\in\mathcal S} \Shadow(\chapeau x,t,R).
\end{align*}
The event $\Shadow(\chapeau x,t, R)$ means that each infinite path starting from $\eta^1(t)$ and contained in the cylinder $\Cyl_+(\chapeau x ,R)$ necessarily meets $\partial\eta(t)\cap\eta^2(t)$.
Loosely speaking, on the event $\Shadow(\chapeau x,t, R)$, the strong infections casts a shadow of radius $R$ on the weak infection. 

Our main result says that if the strong infection occupies a too large portion of the frontier, \ie if $\Shade(t, Rt^{1/2+\eta})$ occurs, then the survival probability of the slow species $1$ is very small:

\begin{theorem}
\label{avantage3}
Consider $M>0$ and $\eta \in (0,1/2)$. There exist two strictly positive constants $A,B$ such that 
$$
\forall t\ge 0 \quad \P \left( \mathcal{G}^1\cap \Shade(t, Mt^{1/2+\eta}) \right)\le A\exp(-Bt^\eta).
$$
\end{theorem}

\begin{coro}
\label{avantage4}
Define $R_t$ as the supremum of the $r$ for which $\Shade(t, r)$ occurs. Let $\eta>0$.
Then, on the event $\mathcal{G}^1$, we almost surely have
$$\lim_{t\to +\infty} \frac{R_t}{t^{1/2+\eta}}=0.$$
\end{coro}

Remember that, by Lemma~\ref{coupling} and the asymptotic shape result for first-passage percolation, the diameter of $\eta(t)$ is of order $t$. 

To obtain the absence of strong coexistence, we also need a control on the number of such stains of species $2$ on the surface of $\eta(t)$: in dimension larger or equal to three, the set $\partial\eta(t)\cap\eta^2(t)$ is not necessarily connected.
On the contrary, in dimension two, the set $\partial\eta(t)\cap\eta^2(t)$ is connected, which enables us to prove the absence of strong coexistence. Consider any norm $\|.\|$ on $\R^2$ and its discrete balls $\mathcal B(t)=\{y \in \Z^2: \; \|y\| \le t\}$:

\begin{theorem} 
\label{th-densite}
For the two dimensional lattice, we have \\ 
1. For every $\beta\in (0,1/2)$, there exists a constant $C>0$ such that, almost surely on the event $\mathcal G^1$
$$ \forall t>0 \quad
\frac{| \eta^2(\infty) \cap \mathcal B(t)|}{|\mathcal B(t)| }  \le  \frac{C}{t^{1/2-\beta}}.
$$
2. For every $\beta>0$, almost surely on the event $\mathcal G^1$
$$\miniop{}{\limsup}{t\to +\infty} \frac{ \text{Diam}\left(\eta^2(\infty)+[-1/2,1/2]^2 \right) \cap \mathcal S_{p_1}(t))}{t^{1/2+\beta}}=0.$$
3. Strong coexistence almost surely does not happen.
\end{theorem}

The next corollary of Theorem~\ref{avantage2} precises Lemma 5.2 in H\"aggstr\"om-Pemantle~\cite{Haggstrom-Pemantle-2}: when coexistence occurs, the two species globally grow with the speed of the slow species, as if the slow species were alone. It also corresponds to a weak version of moderate deviations for first-passage percolation (see the results by Kesten and Alexander, recalled as Proposition~\ref{shapeMD} in the next section).

\begin{theorem}
\label{pasdexcroissance} Let $\beta>0$ and $\eta \in (0,1/2)$. There exist two strictly positive constants $A,B$ such that for every $t \ge 0$:
$$\P \left( \mathcal G^1 \backslash \left\{\mathcal B_{p_1} \left( t - \beta t^{1/2+\eta} \right) \subset \eta(t) \subset \mathcal B_{p_1} \left( t + \beta t^{1/2+\eta} \right) \right\} \right) \le A \exp(-Bt^\eta).$$
\end{theorem}

The estimates we obtained in this paper finally allows us to recover the "except perhaps for a denumerable set" non-coexistence result by H\"aggstr\"om and Pemantle, and to extend it to more general families of passage times indexed by a continuous parameter:
\begin{theorem}
\label{HPgeneral}
Let $(\nu_p)_{p\in I}$  be a family of probability measures indexed by a subset of $\R$. We assume that 
\begin{enumerate}
\item for each $p \in I$, $\nu_p(0) < p_c$
\item for each $p \in I$, $\nu_p(\inf \supp \nu_p) < \pcfleche$
\item for each $p \in I$, there exists $\gamma>0$ such that $\int_{[0, \infty)}\exp(\gamma x)\ d\nu_p(x)<\infty$,
\item for each $p,q \in I$, $p < q \Rightarrow \nu_p \gstoch \nu_q \text{ and } \nu_p \neq \nu_q$,
\item for each $p,q \in I$, $\forall k,l \in \N, \; \nu_{p}^{*k}\otimes \nu_{q}^{*l}(\{(x,x): \, x \in [0,+\infty)\})=0.$
\end{enumerate}
Denote by $\P_{p,q}$ the law of the competition process where species $1$ (\resp $2$) uses passage times with law $\nu_p$ (\resp $\nu_q$).
Then for each fixed $q \in I$, the set $$\{p \in I: \; p\le q \text{ and } \P_{p,q}(\Coex)>0\}$$ is a subset of the points of discontinuity of the non-decreasing map
$p\mapsto \P_{p,q}(\mathcal G^1)$,
and is therefore at most denumerable.
\end{theorem}

\subsection*{Organization of the paper}
The paper is organized as follows: in Section~\ref{preliminaires}, we give a series of useful results in first-passage percolation. Most of them are classical
 and are recalled without proof. We also give there the proof of Lemma~\ref{coupling}.
Section~\ref{core} is mainly devoted to the proof of Theorem~\ref{avantage4}, which is the technical core of the paper.  Section~\ref{sectionmoderate} establishes Theorem~\ref{pasdexcroissance}).
In Section~\ref{sectiondimdeux}, we improve for the two dimensional lattice the results of Section~\ref{core} into the more friendly Theorem~\ref{th-densite}. 
The last section extends the Häggström-Pemantle Theorem to the 
present context as announced in Theorem~\ref{HPgeneral}.

\section{Preliminary results}
\label{preliminaires}
\subsection{First-passage percolation results}

Let us recall some classical results about simple first-passage 
percolation. We assume here that the passage times are independent identically distributed  with common law
$\nu$ satisfying 
\begin{itemize}
\label{hypo}
\item $\nu(0)<p_c$; 
\item for some $\gamma>0$, $\int_{[0, +\infty)}\exp(\gamma x)\ d\nu(x)<+\infty$.
\end{itemize}
Denote by $\norme{.}{\nu}$ the norm given by the shape theorem, and by $\mathcal B^x(t)$ the discrete ball relatively to $\norme{.}{\nu}$   with center $x$ and radius $t$. The first two results give large deviations and moderate deviations for fluctuations with respect to the asymptotic shape, and the third one gives the strict monotonicity result for the asymptotic shape with respect to the distribution of the passage times:

\begin{prop}[Grimmett-Kesten~\cite{grimmett-kesten}] 
\label{shapeGD}
For any $\epsilon>0$, there exist two strictly positive constants $A,B$ such that
$$\forall t>0 \quad \P \left(
\mathcal B^0 ((1-\epsilon)t)\subset B^0(t) \subset \mathcal B^0((1+\epsilon)t)
\right) \ge 1-  A \exp(-B t).$$
\end{prop}

\begin{prop}[Kesten~\cite{kesten-modere},Alexander~\cite{Alex97}]
\label{shapeMD}
 For any $\beta>0$, for any $\eta \in (0,1/2)$, 
there exist two strictly positive constants $A,B$ such that
$$\forall t>0 \quad \P \left(
 \mathcal B^0 (t -\beta t^{1/2+\eta})\subset B^0(t) \subset  \mathcal B^0(t +\beta t^{1/2+\eta})
\right) \ge 1-  A \exp(-B t^\eta).$$
\end{prop}

\begin{prop}[van den Berg-Kesten \cite{vdb-kes}]
\label{strictecomp}
Let $\nu_{p_1}$ and $\nu_{p_2}$ be two probability measures on $[0,+\infty)$ 
satisfying $(H_1)$, $(H_3)$, $(H_4)$, $(H_5)$.

There exists a constant $C_{p_1,p_2}\in (0,1)$ such that 
$$\forall x \in \R^d \quad \norme{x}{p_2} \le C_{p_1,p_2} \norme{x}{p_1}.$$
\end{prop}

Note that in \cite{vdb-kes}, the proof of this result is only written for the time constants. Nevertheless, it applies in any direction and computations
can be followed in order to preserve a uniform control, whatever direction one
considers. See for instance Garet and Marchand~\cite{GM-bernoulli} 
for a detailed proof in an analogous situation. In the same way, the large deviation result Proposition~\ref{shapeGD} is only stated in~\cite{grimmett-kesten} for the time constant, but the result can be extended uniformly in any direction, as it is done in Garet and Marchand~\cite{GM-large} for chemical distance in supercritical Bernoulli percolation.
As far as Proposition~\ref{shapeMD} is concerned, it is a by-product of the proof of Theorem~3.1 in Alexander~\cite{Alex97}. We include here a short proof for convenience.

\begin{proof}[Proof of Proposition~\ref{shapeMD}]
The outer bound for $B^0(t)$ follows from Kesten~\cite{kesten-modere}, Equation~(1.19): there exist positive constants $A_1,B_1$ such that  for all $t>0$, we have
$$\P( B^0(t) \not\subset  \mathcal B^0(t +\beta t^{1/2+\eta})\le A_1\exp(-B_1t^{\eta}).$$
Turning to the inner bound, we follow the lines of Alexander's proof:
for $x,y\in\Zd$, let us define the travel time between $x$ and $y$ by 
$T(x,y)=\inf\{t\ge 0:\ y\in B^x(t)\}.$
Then we have:
\begin{eqnarray*}
&& \P(B^0(t)\not\supset\mathcal B^0(t -\beta t^{1/2+\eta})) \\
& \le &\P(B^0(t)\not\supset\mathcal B^0(t/2))+\miniop{}{\sum}{t/2\le \norme{x}{\nu}\le t-\beta t^{1/2+\eta}}\P(T(0,x)\ge t)\\
 & \le & A\exp(-Bt)+\miniop{}{\sum}{t/2\le \norme{x}{\nu}\le t-\beta t^{1/2+\eta}}\P(T(0,x)\ge t),\\
\end{eqnarray*} 
where $A$ and $B$ are determined by Proposition~\ref{shapeGD} with $\epsilon=1/2$.
By Alexander~\cite{Alex97}, Theorem~3.2, there exist positive constants $C'_4,M$ such that
$$\norme{x}{\nu}\ge M\Longrightarrow \E T(0,x)\le \norme{x}{\nu}+C'_4\norme{x}{\nu}\log \norme{x}{\nu}.$$
Assume $t/2\ge M$ and let $x$ with $t/2\le \norme{x}{\nu}\le t-\beta t^{1/2+\eta}$.
Let $C_M=\inf\{2^{1/2+\eta}\beta-C'_4 y^{-\eta}\log y: \ y\ge M\}$.
We assume that $M$ is so large that $C_M>0$.
We have $\E T(0,x)\le t-C_M\norme{x}{\nu}^{1/2+\eta}\le t-C'_M \|x\|_2^{1/2+\eta}$, where $C_M'$ is a positive constant, and then
$$\P(T(0,x)\ge t)\le \P(T(0,x)-\E T(0,x)\ge C'_M \|x\|^{1/2+\eta}).$$
By Kesten's result~\cite{kesten-modere}, Equation~(2.49) (see also Equation~(3.7) in Alexander~\cite{Alex97}), there exist positive constants $C_5,C_6$ such that
$$ \P(T(0,x)-\E T(0,x)\ge C'_M \|x\|^{1/2+\eta})\le C_5\exp(-C_6\|x\|_2^{\eta}),$$ provided that $M$ is large enough.
Finally, it gives that
$$\miniop{}{\sum}{t/2\le \norme{x}{\nu}\le t-\beta t^{1/2+\eta}}\P(T(0,x)\ge t)
\le \Card{\mathcal B^0(t)} C_5\exp(-C_6 (t/2)^{\eta})\le C_5'\exp(-C_6' t^{\eta}),$$
where $C_5', C_6'$ are positive constants.
\end{proof}

The next lemma ensures that the minimal time needed to cross the cylinder $\Cyl_z(\chapeau x,h,r)$ from bottom to top, using only edges in the cylinder, can not be much larger than the expected value $h\norme{\chapeau x}{\nu}$.

\begin{lemme}
\label{GDlong}
For $z\in \Rd$, $\chapeau x\in \mathcal S$, and $h,r>0$ large enough, we can define the point
 $s_0$ (the point $s_f$) to be the integer point in $\Cyl_z(\chapeau x,r,h)$ which is 
closer to $z$ (respectively,  $z+h\chapeau x$).
We define the crossing time $t[\Cyl_z(\chapeau x,h,r)]$ of the cylinder
$\Cyl_z(\chapeau x,r,h)$ as the minimal time needed to cross it from $s_0$ to $s_f$, using only edges in the cylinder. 

Then for any $\epsilon>0$, and any function $f:\R_+\to \R_+$ with 
$\lim_{+\infty} f=+\infty$, 
there exist two strictly positive constants $A$ and $B$ such that
$$
\begin{array}{l}
\forall z\in\Rd\quad \forall\chapeau x \in \mathcal S \quad \forall h>0 \\
  \P \left(
t[\Cyl_z(\chapeau x,h,f(h))] \ge \norme{\chapeau x}{\nu} ( 1+ \epsilon) h
\right) \le A \exp (-Bh).
\end{array}$$
\end{lemme}

Note that this gives the existence of a nearly optimal path from $z$ to $z+h\chapeau x$ that remains at a distance less than $f(h)$ of the straight line. This result can be interesting on its own as we often miss information on the position of the real optimal paths.

\begin{proof}
For $x,y\in\Zd$, denote by $I_{x,y}$ the length of the shortest path from
$x$ to $y$ which is inside $\mathcal{B}_x(1,25\norme{x-y}{\nu})\cap \mathcal{B}_y(1,25\norme{x-y}{\nu})$. Of course $I_{x,y}$ as the same law that $I_{0,x-y}$, and we simply write $I_x=I_{0,x}$.
We begin with an intermediary lemma.

\begin{lemme}
\label{tfl}
Let $\epsilon,a$ in $(0,1)$ and $\|.\|$ be any norm on $\Rd$. There exists $M_0$ such that 
for each $M\ge M_0$, there exist  $\rho \in (0,1)$ and
$t>0$ such that
$$\|x\|\in [aM,M/a] \Longrightarrow \E \exp(t (I_x-(1+\epsilon)\norme{x}{\nu}))\le \rho.$$
\end{lemme}

\begin{proof}[Proof of Lemma~\ref{tfl}]
Note that by norm equivalence, we can restrict ourselves to $\|.\|_1$. Let $Y$ be a random variable with law $\nu$ and let $\gamma>0$ be such that $\E e^{2\gamma Y}<+\infty$.

First, the large deviations result, Proposition~\ref{shapeGD}, easily implies the following almost sure convergence:
$$\miniop{}{\lim}{\|x\|_1\to +\infty} \frac{I_x}{\norme{x}{\nu}}=1.$$
By considering a deterministic path from $0$ to $x$ with length $\|x\|_1$, we see that $I_x$ is dominated by a sum of $\|x\|_1$ independent copies of $Y$ denoted by $Y_1,\dots, Y_{\|x\|_1}$, and thus
$I_x/{\|x\|_1}$ is dominated by $$\displaystyle \frac1{\|x\|_1}\sum_{k=1}^{\|x\|_1} Y_i.$$
This family is equi-integrable by the law of large numbers. So $\displaystyle (T_x/{\|x\|_1})_{x\in\Zd\backslash\{0\}}$ and  then $\displaystyle (I_x/{\norme{x}{\nu}})_{x\in\Zd\backslash\{0\}}$ are also equi-integrable families, which implies that 
\begin{equation}
\label{equiinteg}
\miniop{}{\lim}{\|x\|_1\to +\infty} \frac{\E I_x}{\norme{x}{\nu}}=1.
\end{equation}

Note now that for every $y \in \R$ and $t\in (0,\gamma]$, 
$$e^{ty} \le  1+ty +\frac{t^2}2y^2 e^{t|y|}
 \le   1+ty+\frac{t^2}{\gamma^2} e^{2\gamma|y|}.$$
Define $\tilde{I}_x=I_x-(1+\epsilon)\norme{x}{\nu}$ and suppose that $t\in (0,\gamma]$. Then, as $|\tilde{I}_x|\le I_x+2 \norme{x}{\nu}$, the previous inequality implies that
$$
e^{t\tilde{I}_x} \le   1+t\tilde{I}_x +\frac{t^2}{\gamma^2} e^{4\gamma\norme{x}{\nu}}e^{2\gamma {I}_x}.
$$
As $\norme{x}{\nu}\le \|x\|_1\norme{e_1}{\nu}$ and $I_x \le Y_1+ \dots +Y_{\|x\|_1}$, if we set $R=e^{4\gamma\norme{e_1}{\nu}}\E e^{2\gamma Y}$, we obtain that
$$\E e^{t\tilde{I}_x}\le 1+t \left( \E \tilde{I}_x+\frac{t}{\gamma^2}R^{\|x\|_1} \right).$$
Considering Equation (\ref{equiinteg}), let $M_0$ be such that
$\|x\| \ge aM_0$ implies $\frac{\E I_x}{\norme{x}{\nu}}\le 1+\epsilon/3$.\\
For $x$ such that $\|x\|_1\ge aM_0$, we have $\E \tilde{I}_x \le -\frac23\epsilon \norme{x}{\nu}$,
so
$$
\E e^{t\tilde{I}_x} \le  1+t \left( -\frac23\epsilon\norme{x}{\nu}+\frac{t}{\gamma^2}R^{\|x\|_1} \right)
 .
$$
Therefore, we take $t=\min \left\{ \gamma,\gamma^2\norme{e_1}{\nu}\frac{\epsilon}3R^{-M/a} \right\}>0$ and $\rho=1-\frac13\epsilon\norme{e_1}{\nu}t<1$.
\end{proof}

Let us come back now to the proof of Lemma~\ref{GDlong}. Let $\epsilon \in (0,1)$ and consider the integer $M_0 \in \mathbb{N}$ given by Lemma~\ref{tfl}. As $\norme{.}{\nu}$ is a norm, there exist two strictly positive constants $c$ and $C$ such that
\begin{equation}
\label{lemm1E}
\forall\chapeau x \in \mathcal S \quad c \le \norme{\chapeau x}{\nu} \le C.
\end{equation}
Let $M_1 \ge M_0$ be an integer large enough to have 
\begin{equation}
\label{lemm1A}
(1+2\epsilon)\ge (1+\epsilon)\left(1+\frac{4\norme{e_1}{\nu}}{cM_1}\right).
\end{equation}

Consider $h>M_1$ and set $N=1+\text{Int}(h/M_1)$ -- $\text{Int} (x)$ denotes the integer part of $x$ -- and, for each $i\in\{0,\dots,N\}$
denote by $x_i$ the integer point in the cylinder which is the closest to
$z+\frac{ih\chapeau x}{N}$. Note that
\begin{equation}
\label{lemm1B}
\forall h \ge M_1 \quad  \frac{M_1}2 \le\left(1-\frac1N\right) M_1 \le \frac h N \le M_1.
\end{equation}
and that for each $i, j\in\{0,\dots,N-1\}$, 
\begin{equation}
\label{lemm1C}
\left| \norme{x_i-x_{j}}{\nu}-|j-i| \frac{h\norme{\chapeau x}{\nu}}N \right|\le 2\norme{e_1}{\nu}.
\end{equation}

1. Applying (\ref{lemm1C}), (\ref{lemm1E}) and (\ref{lemm1B}), we obtain that for each $i \in \{0, \dots, N-1\}$, 
\begin{eqnarray*}
\frac{h\norme{\chapeau x}{\nu}}N -2\norme{e_1}{\nu} \le & \norme{x_i-x_{i+1}}{\nu} & \le \frac{h\norme{\chapeau x}{\nu}}N +2\norme{e_1}{\nu} \\
c\frac{M_1}2 -2\norme{e_1}{\nu} \le & \norme{x_i-x_{i+1}}{\nu} & \le CM_1+2\norme{e_1}{\nu}.
\end{eqnarray*}
So we can find $a>0$ such that, by increasing $M_1$ if necessary, 
\begin{equation}
\label{lemm1D}
\forall i \in \{0, \dots, N-1\} \quad aM_1\le \norme{x_i-x_{i+1}}{\nu} \le M_1/a.
\end{equation}

2. Let $h_1 \ge 0$ be such that
$$\forall h \ge h_1 \quad  f(h)\ge 2,5C(M_1+1)\norme{e_1}{\nu}+1.$$
If we take now $h$ larger than $h_1$, and  if $y \in \mathcal{B}^{x_i}(1,25\norme{x_i-x_{i+1}}{\nu})$ for some $i \in\{0,\dots,N-1\}$, then, with (\ref{lemm1E}) and (\ref{lemm1C}),
\begin{eqnarray*}
d(y-z,\mathbb R\chapeau x) & = & d(y,z+\mathbb R\chapeau x)  \le \left\|y-z-\frac{ih}N\chapeau x \right\|_2 
 \le \|y-x_i\|_2+\left\|x_i-z-\frac{ih}N\chapeau x \right\|_2 \\
& \le & 1,25C\norme{x_i-x_{i+1}}{\nu} +1 
 \le  1,25C \left( \frac{h\norme{\chapeau x}{\nu}}N +2\norme{e_1}{\nu}\right) +1.
\end{eqnarray*}
But $h/N \le M_1$ and $\norme{\chapeau x}{\nu} \le \|\chapeau x \|_1 \norme{e_1}{\nu} \le 2\norme{e_1}{\nu}$, thus
$$d(y-z,\mathbb R\chapeau x) \le f(h).$$
On the other hand, 
\begin{eqnarray*}
\langle y-z,\chapeau x \rangle & = & \langle y-x_i,\chapeau x \rangle + \langle x_i-\left( z+ \frac{ih}N\chapeau x \right),\chapeau x \rangle + \langle \frac{ih}N\chapeau x,\chapeau x \rangle \\ 
\text{\ie } \left| \langle y-z,\chapeau x \rangle - \frac{ih}N \right| & \le & \left\|  y-x_i \right\|_2 +1 \le 2,5C(M_1+1)\norme{e_1}{\nu}+1.
\end{eqnarray*}
We choose then $i_0 \in \mathbb N$ such that:
$$i_0 \ge \frac{2}{M_1}(2,5C(M_1+1)\norme{e_1}{\nu}+1).$$
Then, if $h_2$  is such that $1+\text{Int}(h_2/M_1)\ge 3i_0$, we obtain:
\begin{equation}
\label{lemm1H}
\forall h \ge h_2 \quad \forall i \in \{i_0, \dots, N-1-i_0\} \quad \mathcal{B}^{x_i}(1,25\norme{x_i-x_{i+1}}{\nu}) \subset \Cyl_z(\chapeau x,R,h).
\end{equation}

3. There exists a deterministic path inside the cylinder from $x_0$ to $x_{i_0}$
(from $x_{N-i_0}$ to $x_N$) which uses less than $i_0 \frac{h\|\chapeau x\|_1}N +2$
edges: we denote by
$L_{start}$ (respectively,   $L_{end}$) the random length of this path. By Equation (\ref{lemm1B}), we have
$$\forall h \ge M_1 \quad i_0 \frac{h\|\chapeau x\|_1}N +2 \le i_0 \frac{2h}N +2 \le 2(i_0+1)M_1.$$
If $h>h_3=\frac{3(i_0+1)M_1 \E Y }{\epsilon\norme{\chapeau x}{\nu}}$, Chernoff's theorem gives the existence of two strictly positive constants $A_1,B_1$ such that  
\begin{equation}
\label{lemm1G}
\forall h>0 \quad \P \left(L_{start}>{\epsilon} h \norme{\chapeau x}{\nu} \right)+\P\left(L_{end}>{\epsilon} h\norme{\chapeau x}{\nu} \right)\le A_1e^{-B_1 h}.
\end{equation}

4. So, provided that $h\ge h_2$, we have by (\ref{lemm1G}), inside the cylinder, a path from $x_0$ to $x_N$
with length
$$L_{start}+\sum_{i=i_0}^{N-i_0-1} I_{x_i,x_{i+1}}+L_{end}.$$
By Equation (\ref{lemm1C}), if $h$ is large enough, for each $i, j\in\{0,\dots,N-1\},$ 
\begin{align*}
& \mathcal{B}^{x_i}(1,25\norme{x_i-x_{i+1}}{\nu}) \cap \mathcal{B}^{x_{i+1}}(1,25\norme{x_{i}-x_{i+1}}{\nu}) \\
&\cap \mathcal{B}^{x_j}(1,25\norme{x_j-x_{j+1}}{\nu})\cap \mathcal{B}^{x_{j+1}}(1,25\norme{x_{j}-x_{j+1}}{\nu})=\varnothing
\end{align*}
as soon as $|j-i|\ge 2$.
We thus introduce, for $j \in \{0,1\}$, the sums:
$$S_{j}=\sum_{\substack{I\le i\le N-I-1   \\ i=j \mod 2}}  I_{x_i,x_{i+1}}.$$
Note that, with (\ref{lemm1C}), (\ref{lemm1E}) and  (\ref{lemm1A}) for each $j \in \{0,1\}$,
\begin{eqnarray*}
\sum_{\substack{ i_0\le i\le N-i_0-1  \\ i =j \mod 2}}    \norme{x_{i+1}-x_i}{\nu} 
& \le & \frac{N-2I}2 \left( \frac{h \norme{\chapeau x}{\nu}}N +2 \norme{e_1}{\nu} \right) 
 \le  \frac{h \norme{\chapeau x}{\nu}}2 \left( 1 + \frac{2N\norme{e_1}{\nu}}{h \norme{\chapeau x}{\nu}} \right) \\
& \le &  \frac{h \norme{\chapeau x}{\nu}}2 \left( 1 + \frac{4\norme{e_1}{\nu}}{cM_1} \right)
\le \left( \frac{1+2\epsilon}{1+\epsilon} \right) \frac{h \norme{\chapeau x}{\nu}}2.
\end{eqnarray*}
Then, by independence of the terms in each $S_j$,
\begin{eqnarray*}
\P \left( S_j\ge \frac{h\norme{\chapeau x}{\nu}}{2}(1+2\epsilon) \right) 
& \le & \P \left( S_j\ge (1+\epsilon)\sum_{\substack{i_0\le i\le N-i_0-1  \\ i =j \mod 2}}    \norme{x_{i+1}-x_{i}}{\nu} \right)\\
& \le & \E \exp \left( t \sum_{\substack{i_0\le i\le N-i_0-1  \\ i =j \mod 2}}  \left(   I_{x_i,x_{i+1}}-(1+\epsilon)\norme{x_{i+1}-x_{i}}{\nu} \right)\right)\\
& \le &\prod_{\substack{i_0\le i\le N-i_0-1  \\ i =j \mod 2}} \E \exp(t(I_{x_i,x_{i+1}}-(1+\epsilon)\norme{x_{i+1}-x_{i}}{\nu}))\\
\end{eqnarray*}
By (\ref{lemm1D}), for each $i$, we have $\norme{x_i-x_{i+1}}{\nu}\in [aM_1,M_1/a]$, so we can apply Lemma~\ref{tfl}: there exists some $\rho<1$, such that for every $h$ large enough,
\begin{equation*}
\label{lemm1F}
\forall j \in \{0,1\} \quad \P\left( S_j\ge \frac{h \norme{\chapeau x}{\nu}}2 (1+2 \epsilon) \right)  \le  \rho^{N/2} \le  \rho^{h/(2M_1)}.
\end{equation*}
Together with (\ref{lemm1G}), this proves the estimate of the lemma.
\end{proof}

\subsection{Comparisons with first-passage percolation}

We now prove Lemma~\ref{coupling}, using an algorithmic building analogous to the one used to define the competition process in the introduction.

\begin{proof}[Proof of Lemma~\ref{coupling}]
The inclusion $\eta^1(t) \subset B_{p_1}^{s_1}(t)$ is obvious: by construction of the process, if $x \in \eta^1(t)$, there exists a path between $s_1$ and $x$ included in $\eta^1(t)$, and whose travel time  is thus less than $t$. The second inclusion $\eta^2(t)  \subset  B_{p_1}^{s_2}(t)$ is proved in the same way.

Let us now prove the third inclusion $B_{p_1}^{s_1}(t) \subset \eta(t).$
We are going to build the first-passage process by Dijkstra's algorithm, in a formalism analogous to the one used to define the competition process.

Recall that $\Omega=([0,+\infty)\times[0,+\infty))^{\Ed}$ is endowed with the measure $\P=m^{\otimes\Ed}$. Consider a fixed configuration $\omega \in \Omega$.

Let  $E'=[0,+\infty]^{\Zd}$.
We recursively define a $E'$-valued sequence $(X'_n)_{n\ge 0}$ and a non-negative
sequence  $(T'_n)_{n\ge 0}$. The sequence $(T'_n){n\ge 0}$ contains the successive times of infections, while a point $\epsilon=(\epsilon(z))_{z \in \Zd} \in E$ codes, for each site $z$, its times of infection $\epsilon(z)$. We start the process with the single source $s_1$, and set $T'_0=0$ and 
$$
X'_0=(X_0(z))_{z \in \Zd} \text{ with } X'_0(s_1)=0 \text{ and } X'_0(z)=+\infty \text{ if } z \neq s_1.
$$
Then, for $n \ge 0$, define the next time of infection:
$$T'_{n+1}=\inf\{ X'_n(y)+\omega_{\{y,z\}}^1: \quad \{y,z\}\in\Ed\}.$$
The infimum is reached for some couples $(y_i,z_i)$, meaning that the $z_i$ are being infected from the $y_i$:
$$\forall x \in \Zd\backslash \{z_i,i\}, \; X'_{n+1}(x)=X'_n(x) \text{ and } X'_{n+1}(z_i)=X'_n(y_i)+\omega_{\{y_i,z_i\}}^1.$$
We also note that $\eta'(t)$, the set of infected points at time $t$ by 
$$\forall n \in \N \quad \forall t\in [T_n,T_{n+1})\quad  \eta'(t)=\{z\in\Zd: \; X'_n(z)<+\infty\}.$$
By Dijkstra's algorithm, $\eta'(t)$ is exactly the set $B_{p_1}^{s_1}(t)$.

We now proceed by induction to prove that for every $n \in \N$ 
$$ (H_n) \quad \forall x \in \Zd \quad X'_n(x)\ge \min(X^1_n(x), X^2_n(x)).$$
Clearly, $(H_0)$ is true. Assume that $(H_n)$ holds. We have the following alternative:
\begin{itemize}
\item If $X'_{n+1}(x)=+\infty$, it is obvious that $X'_{n+1}(x)\ge \min(X^1_{n+1}(x), X^2_{n+1}(x)).$
\item If $X'_{n}(x)<+\infty$ then $X'_{n+1}(x)=X'_{n}(x)$ and, as $X'_n(x)\ge \min(X^1_n(x), X^2_n(x))$, the number $\min(X^1_n(x), X^2_n(x))$ is also finite, and thus $\min(X^1_n(x), X^2_n(x))=\min(X^1_{n+1}(x), X^2_{n+1}(x))$ -- recall that $X^1_n(x)$ and $X^2_n(x)$ can not be simultaneously finite. Consequently, $X'_{n+1}(x)\ge \min(X^1_{n+1}(x), X^2_{n+1}(x)).$
\item If $X'_{n+1}(x)<+\infty$ and $X'_{n}(x)=+\infty$, the point $x$ is being infected at time $T'_{n+1}$  through the edge $e$ from the point $y$, which is consequently such that $X'_n(y)<+\infty$. As $X'_n(y)\ge \min(X^1_n(y), X^2_n(y))$, and $\omega^1_e \ge \omega^2_e$, 
$$X'_{n+1}(x)=X'_n(y)+\omega^1_e \ge \min(X^1_{n+1}(y), X^2_{n+1}(y)).$$
\end{itemize}  

Note that
\begin{eqnarray*}
\eta (t) & = & \{ z \in \Zd: \; \exists n \in \N, \; \min(X^1_n(z),X^2_n(z)) \le t \}, \\
\eta'(t) & = & \{ z \in \Zd: \; \exists n \in \N, \; X'_n(z) \le t \}.
\end{eqnarray*}
It is then obvious that $\eta'(t) \subset \eta (t)$.
\end{proof}

\section{Coexistence can not be observed by a medium resolution satellite}
\label{core}

The proof of Theorem~\ref{avantage3} follows, at least in its main lines,
the strategy initiated by 
Häggström and Pemantle: the aim is to prove
that an event, suspected  to be incompatible with the survival of the weak,
 allows  the strong to  grant themselves, with high probability, a family of shells that surround the weak, preventing thus coexistence.
In the Häggström-Pemantle paper~\cite{Haggstrom-Pemantle-2}, the objective is to show that when coexistence occurs, it is unlikely that the strong can advance
significantly beyond the weak.
The event considered here    is of a different nature: we must prove that the strong can not occupy a too large region on the frontier of the infected zone. Obviously, this requires finer controls. Moreover, the use of the shape theorem is
 not sufficient: moderate deviations for
the fluctuations with respect to the asymptotic shape provide sharper estimates. Some more technical difficulties also
follow from the loss of some nice properties of exponential laws.
However, this last kind of difficulties has already been overcome by the authors of 
this paper in the previous article~\cite{GM-bernoulli}. We refer the reader
to this paper for some more comments. 

To prove Theorem~\ref{avantage3}, we first prove an analogous result, Lemma~\ref{avantage2}, in a fixed given direction. Theorem~\ref{avantage4} follows then from a Borel-Cantelli type of argument.

\subsection*{Definitions.}
Let $\mathcal S_{p_2}$ be the unit sphere for the norm $\norme{.}{p_2}$. We define
the \emph{shells}:
for each $A \subset \mathcal{S}_{p_2}$, and every $0<r<r'$, we set
\begin{eqnarray*}
\shell(A,r,r') & = & \{x \in \Zd: \;  x/ \norme{x}{p_2}\in A \text{ and } r \le \norme{x}{p_2}
\le r'\}.
\end{eqnarray*}
So roughly speaking, $A$ is to think about as the set of possible
directions for the points in the shell, while $[r,r']$ is the set of radii.

For $A \subset \mathcal{S}_{p_2}$ and $\phi>0$, define the following
enlargement of $A$:
\begin{eqnarray*}
A \oplus \phi & = & (A+\mathcal{B}^0_{p_2}(\phi))\cap \mathcal{S}_{p_2}.
\end{eqnarray*}

Let us state first three geometric lemmas:
\begin{lemme}
\label{direction}
For any norm $|.|$ on $\Rd$, one has
$$
\forall x,y \in \Rd \backslash\{0\} \quad 
\left| \frac{x}{|x|}-\frac{y}{|y|} \right|  
\le\frac{2|x-y|}{\max\{|x|,|y|\}}.
$$  
\end{lemme}

\begin{lemme}
\label{familledeboules}
For any norm $|.|$ on $\Rd$, there exist a constant $C>0$ such that the unit
sphere for $|.|$ can be covered with $C(1+\frac1{\epsilon})^{d-1}$ balls
of radius $\epsilon$ having their centers on the unit sphere.
\end{lemme}
\begin{proof}
When $|.|=\|.\|_{\infty}$, it is easy to see that the sphere can be covered
with $2d(1+\frac2{\epsilon})^{d-1}$ balls of radius $\epsilon$.

Now let $A,B$ be two strictly positive constants such that
$$\forall x\in\Rd \quad A\|x\|_{\infty}\le |x|\le  B\|x\|_{\infty}$$ and let
 $K=\frac{A}{2B}$.
Suppose, by the previous step, that the unit sphere $\{x\in\Rd:\|x\|_{\infty}=1\}$ is covered by the family of balls $(\mathcal{B}^{x_i}_{\infty}(K\epsilon))_{1 \le i \le n}$ with $\|x_i\|_{\infty}=1$ for each $i$ and $n \le 2d(1+\frac2{K\epsilon})^{d-1}$.
We note $\Psi(x)=\frac{x}{|x|}$.

Let $y$ in the unit sphere for $|.|$: there exists $i \in \{1, \dots,n\}$ such that
$\|x_i-\frac{y}{\|y\|_{\infty}}\|_{\infty}\le K\epsilon$.
Since $\Psi(\frac{y}{\|y\|_{\infty}})=y$, Lemma~\ref{direction} ensures that
$$| \Psi(x_i)-y|\le\frac{2|x_i-\frac{y}{\|y\|_{\infty}}|}{|x_i|}\le\frac{2BK\epsilon}{A}=\epsilon.$$
So, the unit sphere for  $|.|$ can be covered with
 $2d(1+\frac{4B}{A}\frac1\epsilon)^{d-1}$ balls of radius $\epsilon$
having their centers on the unit sphere.
\end{proof}

\begin{lemme}
\label{petitcyl}
Let $R\ge 0$. If $\chapeau u,\chapeau v \in \mathcal S$ are such that
$\|\chapeau u -\chapeau v \|_2\le (\frac{h}{2R})^2$, then
$$\Cyl(\chapeau v,h/2)\cap \mathcal{B}_2(R)\subset \Cyl(\chapeau u,h).$$
Moreover, if $\|\chapeau u -\chapeau v\|_2\le 1/2$, then
$$\Cyl_+(\chapeau v,h/2)\cap \mathcal{B}_2(R)\cap \mathcal{B}_2(h)^c\subset \Cyl_+(\chapeau u,h).$$
\end{lemme}

\begin{proof}
We 
denote $\theta=\|\chapeau u-\chapeau v\|_2$.
Then $d(y,\R\chapeau u)^2=\| y-\langle y,\chapeau u\rangle\chapeau u\|_2^2=\|y\|_2^2-\langle y,\chapeau u\rangle^2$
and 
\begin{eqnarray*}
d(y,\R\chapeau u)^2- d(y,\R\chapeau v)^2& = & \langle y,\chapeau v\rangle^2-\langle y,\chapeau u\rangle^2\\ 
& \le & \langle y,\chapeau v-\chapeau u\rangle \langle y,\chapeau u+\chapeau v\rangle
  \le  2\theta \|y\|_2^2.
\end{eqnarray*}
 
Suppose first that $y\in \Cyl(\chapeau v,h/2)\cap \mathcal{B}_2(R)$: we have
$$d(y,\R\chapeau u)^2- d(y,\R\chapeau v)^2\le 2\theta R^2\le h^2/2.$$
So $ d(y,\R\chapeau u)^2=d(y,\R\chapeau v)^2 +( d(y,\R\chapeau u)^2- d(y,\R\chapeau v)^2)\le h^2/4+h^2/2 \le h^2$, which means that $y\in  \Cyl(\chapeau u,h)$.

Now suppose that $\theta\le 1/2$ and $\|y\|_2\ge h$. We have
$d_2(y,\R\chapeau v)=\|y-\langle y,v\rangle\chapeau v\|_2\le h/2$ and
$\langle y,\chapeau v\rangle\ge 0$, so $\langle y,\chapeau v\rangle=|\langle y,\chapeau v\rangle|\ge \|y\|_2- h/2$. This implies
$$\langle y,\chapeau u\rangle=
\langle y,\chapeau v\rangle -\langle y,\chapeau v-\chapeau u\rangle
\ge \langle y,\chapeau v\rangle-\|y\|_2\theta\ge \|y\|_2(1-\theta)-h/2\ge \|y\|_2\left(\frac12-\theta\right)\ge 0,$$
which ends the proof.
\end{proof}

The next lemma ensures that if $\Shadow(\chapeau x, t, Rt^{1/2+\eta})$, then the strong infection manages with high probability to colonize a small shell which gives it a positional advantage.
Let $K_1>0$ be such that 
\begin{equation}
\label{K1}
\forall x \in \Rd \quad \|x\|_2 \le K_1 \norme{x}{p_1}.
\end{equation}
\begin{lemme}
\label{petit-bout-de-coquille}  
Let $R>0$, $\eta \in (0,1/2)$ and  $\delta \in (0, R/K_1)$.
Choose then $\delta'>0$ with $\delta< \delta'<\min \left\{ \delta/C_{p_1,p_2}, R/K_1 \right\}$.
Choose  $\gamma,\gamma'$ such that $1<\gamma<\gamma'$ and 
$\theta>0$.
For any $\chapeau x \in \mathcal S$, any $t>0$, any $\gamma,\gamma'$ such that $1<\gamma<\gamma'$ and any
$\theta>0$, we define the following events, depending on $\chapeau x, t, \gamma,
\gamma'$ and $\theta$:
\begin{eqnarray*}
E^1& = & \{\eta^1(t+\delta t^{1/2+\eta})\, \subset \, \mathcal B_{p_1}^0( t+\delta't^{1/2+\eta}) \}, \\
E^2& = & \left\{  \eta^2(t+\delta t^{1/2+\eta}) \, \supset \, 
 \shell \left( \left\{ \frac{\chapeau x }{\norme{\chapeau x}{p_2}} \right\} \oplus(\theta t^{-1/2+\eta}),   \frac{\norme{\chapeau x}{p_2}}{\norme{\chapeau x}{p_1}} (t+\gamma \delta' t^{1/2+\eta}), \frac{\norme{\chapeau x}{p_2}}{\norme{\chapeau x}{p_1}} (t+\gamma' \delta' t^{1/2+\eta})  \right) \right\}, \\
 E & = & E^1 \cap E^2.
\end{eqnarray*}
Then there exist $\gamma'_0>1$ and  $\theta_0>0$  such that for any $\gamma, \gamma',\theta$ satisfying 
$1<\gamma<\gamma'<\gamma'_0$ and  $0<\theta<\theta_0$, there exist two strictly positive constants $A,B$ such that 
$$ 
\forall\chapeau x \in \mathcal S \quad \forall t >0 \quad 
\P( \Shadow(\chapeau x, t, Rt^{1/2+\eta}) 
\backslash E ) \le  
A\exp(-Bt^\eta ).  
$$ 
\end{lemme} 
\begin{proof}
Let $R>0$, $\eta \in (0,1/2)$ and  $\delta \in (0, R/K_1)$.
Choose then $\delta'>0$ with  $\delta<\delta'<\min \left\{ {\delta}/{C_{p_1,p_2}}, {R}/{K_1} \right\}$.

\noindent
\underline{Choice of constants.} Choose now $\delta_1, \epsilon$ and $\beta$ positive such that
\begin{align}
\delta' & <\delta_1 < \frac{\delta}{C_{p_1,p_2}}, \label{Adelta1}\\
\delta' &<\delta_1(1-\epsilon), \label{Aepsilon}\\
\beta  & < \min \{\delta'-\delta, \delta_1\epsilon\} . \label{Abeta}
\end{align}
\vspace{0.2cm}
\noindent
\underline{Step 1. Control of the slow $p_1$-infection.} As $\eta<1/2$, 
we have
$$(t +\delta t^{1/2+\eta})+ \beta (t +\delta t^{1/2+\eta})^{1/2+\eta} = t+(\delta+\beta)t^{1/2+\eta} +o(t^{1/2+\eta}).$$
By condition~(\ref{Abeta}), $\delta+\beta<\delta'$. Thus, using the moderate deviation result (Proposition~\ref{shapeMD}), there exist two strictly positive constants $A_1, B_1$ such that $\forall t>0$
\begin{equation}
\P( (E^1)^c) \le A_1 \exp(-B_1 t^\eta). \label{A1} 
\end{equation}
We can thus assume in the following that $E^1$ occurs.

\vspace{0.2cm}
\noindent
\underline{Step 2. Control of the competition process.}
Denote:
\begin{align*}
F_1=F_1(t) =  \left\{ \eta^1(t) \subset  \mathcal B_{p_1}^0(t+\beta t^{1/2+\eta})
\right\} 
\cap  \left\{  \mathcal B_{p_1}^0(t-\beta t^{1/2+\eta}) \subset \eta^1(t) \cup \eta^2(t)
\right\}.
\end{align*}
By Proposition~\ref{shapeMD} and Lemma~\ref{coupling}, there exist two strictly positive constants $A_2,B_2$ such that $\forall t>0$
\begin{equation}
\P( F_1^c) \le A_2 \exp(-B_2 t^\eta). \label{A3} 
\end{equation}
We can thus assume in the following that $F_1$ occurs.

Define an integer approximation of the line $\R\chapeau x$ as follows:
$$D_{\chapeau x}=\{y \in \Zd: \; \exists z  \in \R\chapeau x, \; \|y-z\|_{\infty} \le 1/2\}.$$
Note that $D_{\chapeau x}$ is connected.
Let now $s_0 \in \eta^2(t) \cap D_{\chapeau x}$ be a point which realizes the maximum
$\max_{y \in \eta^2(t) \cap D_{\chapeau x}} \langle y,\chapeau x \rangle.$
On the event $\Shadow(\chapeau x, t, Rt^{1/2+\eta})  \cap F_1 $, the point $s_0$ is well defined and satisfies 
\begin{equation}
\label{A2}
\norme{s_0}{p_1} \ge t-\beta t^{1/2+\eta}.
\end{equation}

\vspace{0.2cm}
\noindent
\underline{Step 3. The weak infection can not fill the hole.}
Define
\begin{align*}
F_2  & =  F_2(t)  =  \left\{ 
\forall x \in \mathcal B_{p_1}(t+\beta t^{1/2+\eta}) \quad B_{p_1}^x(\delta t^{1/2+\eta}) \subset \mathcal B_{p_1}^x(\delta' t^{1/2+\eta})
\right\}, \\
r & =  R-K_1\delta'>0.
\end{align*}
As $\delta'>\delta$, by the large deviation result (Proposition~\ref{shapeGD}), there exist two strictly positive constants $A_3,B_3$ such that 
$\forall t>0$
\begin{equation}
\P(F_2^c)\le \Card{ \mathcal B_{p_1}(t+\beta t^{1/2+\eta})} A_3 \exp(-B_3t^{1/2+\eta})\le A'_3 \exp(-B'_3t^{1/2+\eta}). \label{A4}
\end{equation}
We can thus assume in the following that $F_2$ occurs.
Note that as $\mathcal B_{p_1}^x(\delta') \subset \mathcal B_2^x(K_1\delta')$ by the very definition of $K_1$, the event $\Shadow(\chapeau x, t, Rt^{1/2+\eta})  \cap F_1 \cap F_2$ prevents the weak infection to bother the strong one in its progression after time $t$ inside $\Cyl(\chapeau x, rt^{1/2+\eta})$. 

\begin{figure}[h!]
\scalebox{0.5}{\input{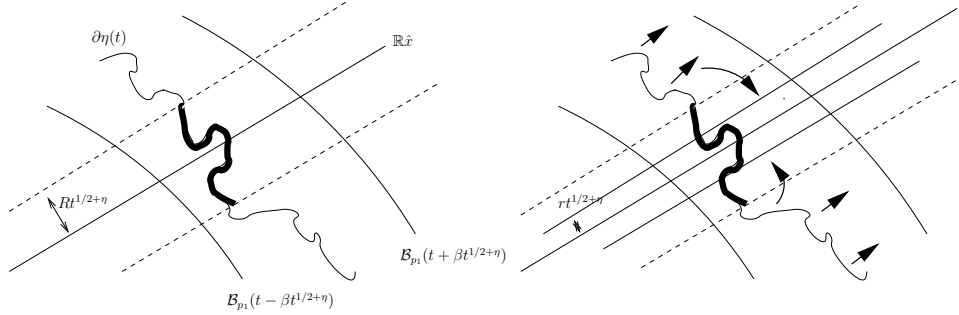}}
\caption{On the left hand side, $\Shadow({x},t, Rt^{1/2+\eta})$: the strong are in bold. On the right hand side, the evolution of the weak: to enter $\Cyl({x}, rt^{1/2+\eta})$, it has to cross a gap of order $(R-r)t^{1/2+\eta}$.} 
\end{figure}

\vspace{0.2cm}
\noindent
\underline{Step 4. The strong manage to escape.}
Remember that, by the choice~(\ref{Aepsilon}), $\delta'<\delta_1(1-\epsilon)$. 
Define $s_f$ as a point in $ D_{\chapeau x } \cap \mathcal B_{p_1}(t+\delta_1(1-\epsilon)t^{1/2+\eta})$ such that $ (s_f +[-1/2,1/2]^d )\cap \mathcal B_{p_1}(t+\delta_1(1-\epsilon)t^{1/2+\eta})^c \neq \varnothing$. Then, by Estimate~(\ref{A2}), we have
\begin{eqnarray*}
\| s_f -s_0\|_2 & \le & \| s_f - \frac{\norme{s_f}{p_1}}{\norme{\chapeau x}{p_1}}\chapeau x \|_2  + \frac{\norme{s_f}{p_1}}{\norme{\chapeau x}{p_1}} - \frac{\norme{s_0}{p_1}}{\norme{\chapeau x}{p_1}}  + \| \frac{\norme{s_0}{p_1}}{\norme{\chapeau x}{p_1}}\chapeau x-s_0\|_2 \\
& \le & 2\sqrt d + \frac{\delta_1(1-\epsilon)+\beta}{\norme{\chapeau x}{p_1}}t^{1/2+\eta} \sim \frac{\delta_1(1-\epsilon)+\beta}{\norme{\chapeau x}{p_1}}t^{1/2+\eta}.
\end{eqnarray*}
Denote
$$F_3=F_3(\chapeau x, t) = \{ s_f \in  \eta^2(t+\delta_1 C_{p_1,p_2}  t^{1/2+\eta})\}.$$
With the definition of $C_{p_1,p_2}$ in Proposition~\ref{strictecomp} and Condition~(\ref{Abeta}), $(\delta_1(1-\epsilon)+\beta)\frac{\norme{\chapeau x}{p_2}}{\norme{\chapeau x}{p_1}} \le C_{p_1,p_2} \delta_1$, and thus 
Lemma~\ref{GDlong} gives two strictly positive constants $A_4,B_4$ such that $\forall\chapeau x \in \mathcal S \quad \forall t >0$
\begin{equation}
 \P \left( \Shadow(\chapeau x, t, Rt^{1/2+\eta})  \cap F_1 \cap F_2 \cap  F_3^c
\right)
 \le  A_4 \exp(-B_4t^{1/2+\eta}) .\label{A6}
\end{equation}

\vspace{0.2cm}
\noindent
\underline{Step 5. The strong colonizes a small shell.} Let $K_2$ be a strictly positive constant such that $\forall x \in \Rd \quad \norme{x}{p_2} \le K_2 \|x\|_2$, and set
\begin{equation}
\label{Adelta2}
\delta_2=\min \left\{ \delta -C_{p_1,p_2} \delta_1, \frac{r}{K_2} \right\}. \nonumber
\end{equation}
Let $\epsilon'>0$ be such that
$ (1+\epsilon')K \delta_2< r$. Denote
$$F_4=F_4(\chapeau x, t)=  \{ 
 \mathcal B_{p_2}^{s_f}( (1-\epsilon')\delta_2 t^{1/2+\eta})\subset \eta_2 (t+\delta t^{1/2+\eta}) \}.$$
Here, the choice we made for $\epsilon'$ ensures that, for $t$ large enough
 $$ \mathcal B_{p_2}^{s_f} ((1+\epsilon')\delta_2 t^{1/2+\eta})\subset ( \mathcal \Cyl(\chapeau x, r t^{1/2+\eta})\cap \mathcal B_{p_1}(t+\delta' t^{1/2+\eta})^c).$$
 Thus, by the large deviation result (Proposition~\ref{shapeGD}), there exist two strictly positive constants $A_5,B_5$ such that
$\forall\chapeau x \in \mathcal S \quad  \forall t>0$, 
\begin{equation}
\P(\Shadow(\chapeau x, t, Rt^{1/2+\eta})  \cap F_1 \cap F_2 \cap F_3 \cap  F_4^c) \le A_5 \exp(-B_5 t^{1/2+\eta}).\label{A7}
\end{equation}
Choose now $\theta_0>0$ and $\gamma'_0>1$ such that
$$
\shell \left( \left\{ \frac{\chapeau x }{\norme{\chapeau x}{p_2}} \right\} \oplus\theta_0t^{-1/2+\eta},   \frac{\norme{\chapeau x}{p_2}}{\norme{\chapeau x}{p_1}} (t+ \delta' t^{1/2+\eta}), \frac{\norme{\chapeau x}{p_2}}{\norme{\chapeau x}{p_1}} (t+\gamma'_0 \delta' t^{1/2+\eta})  \right) \subset
 \mathcal B_{p_2}^{s_f}( (1-\epsilon')\delta_2 t^{1/2+\eta}).
$$
Then $F_4\backslash E^2=\varnothing$, and collecting estimates~(\ref{A1}), (\ref{A3}), (\ref{A4}), (\ref{A6}), and (\ref{A7}), we get the estimate of the lemma.
\end{proof}

The next lemma describes the typical progression of the strong infection from one shell to the next one.

\begin{lemme} 
\label{progressionter}
Let $\phi \in (0,2]$, $h \in (0,1/2)$ and $\alpha \in (1,2)$ be fixed parameters such that
\begin{equation}
(1+h)\left(1+\frac{3\phi}h\right)<\alpha \label{C11}.
\end{equation}
For any $S$ subset of $\mathcal{S}$ and for any $r,s>0$, we define
the following event $E=E(S,r,s)$:
"Any point in the big $\shell(S \oplus \frac{\varphi r}{2(s+r)}\oplus \frac{\varphi (1+h)r}{2(s+(1+h)r)}, s+(1+h)r, s+(1+h)^2 r)$ is linked to a point in the small 
$\shell(S\oplus \frac{\varphi r}{2(s+r)}, s+r, s+(1+h)r)$ by an open path whose length is
less than $\alpha h r$.'' 

Then there exist two strictly positive constants $A$ and $B$, only depending
on $\phi, h, \alpha$, such that for any $r,s>0$ and 
any $S$ of $\mathcal{S}$, 
we have $$\P(E^c) \le A (s+r)^d \exp(-Br).$$
Moreover, we can assume  that all the infection paths needed in $E$ are
completely  
included in the bigger $\shell(S\oplus\frac{(\varphi + 2\alpha h)r}{s+r}, s+[1- 3\phi] (1+h)r, \infty)$.
\end{lemme}

\begin{proof}
Let $\phi \in (0,2]$, $h \in (0,1)$ and $\alpha \in (1,2)$ be fixed parameters satisfying
Equation~(\ref{C11}) and choose, in this order, $\alpha'>1$,
$\epsilon>0$ and $\rho>0$ such that
\begin{align}
& (1+h)\phi+(1+h)^2 (1+ \phi) - (1+h-2\rho) \le \alpha' h < \alpha h, \label{C21} \\
& 4\rho \alpha' h  < \phi \text{ and } 4\rho \alpha' h<\phi(1+h-2\rho), \label{C22} \\
& h-2\rho-\rho \alpha' h  >  0 ,\label{C24} \\
& (1+\epsilon)^2(1+\rho)\alpha' \le  \alpha.\label{C25}
\end{align}
Note that (\ref{C21}) is allowed by (\ref{C11}).
Let $S$ be any subset of $\mathcal{S}$. Denote $T=S \oplus \frac{\varphi r}{2(s+r)}\oplus \frac{\varphi (1+h)r}{2(s+(1+h)r)}$.
Let $ z \in T$: 
$$
\exists v \in S, \; u_1 \in \mathcal{B}_{p_2}^0 \left(\frac{\varphi (1+h)r}{2(s+(1+h)r)}\right), \; u_2 \in \mathcal{B}_{p_2}^0 \left(\frac{\varphi r}{2(s+r)} \right) \text{ such that } z=v+u_1+u_2.
$$
\begin{itemize}
\item As $v \in S$, we have  $\displaystyle v \oplus \frac{\varphi r}{2(s+r)} \subset  S \oplus \frac{\varphi r}{2(s+r)}$.
\item Moreover, $\displaystyle \norme{z-v}{p_2} =\norme{u_1+u_2}{p_2} \le \norme{u_1}{p_2}+\norme{u_2}{p_2} \le \frac{\varphi (1+h)r}{(s+(1+h)r)}$.
\end{itemize}
Thus for any $z \in \shell(T, s+(1+h)r, s+(1+h)^2 r)$, we can choose $\chapeau v_z \in {S}$ such that
\begin{equation}
\norme{\frac{z}{\norme{z}{p_2}} -\chapeau v_z}{p_2} \le \frac{\phi (1+h)r}{s+(1+h)r} \text{ and }\chapeau v_z \oplus \frac{\varphi r}{2(s+r)} \subset  S \oplus \frac{\varphi r}{2(s+r)}.
\nonumber
\end{equation}
We set $v_z=[s+(1+h-2\rho)r]\chapeau v_z$. 
Let us first estimate $\norme{z-v_z}{p_2}$: on the one hand,
\begin{eqnarray}
\norme{z-v_z}{p_2} & \le & \norme{z-\norme{z}{p_2}\chapeau{v}_z}{p_2}+| \, \norme{z}{p_2}- s-(1+h-2\rho)r| \nonumber \\
& \le & \norme{z}{p_2} \frac{\phi (1+h)r}{s+(1+h)r} + \norme{z}{p_2}- s-(1+h-2\rho)r \nonumber \\
& \le &  \norme{z}{p_2} \left(1+ \frac{\phi (1+h)r}{s+(1+h)r} \right)-s-(1+h-2\rho)r \label{modulezvMaj} \\
& \le & s\frac{\phi (1+h)r}{s+(1+h)r}  +[(1+h)^2 \left(1+ \frac{\phi (1+h)r}{s+(1+h)r} \right) - (1+h-2\rho) ] r  \nonumber \\
& \le & [(1+h)\phi+ (1+h)^2 (1+ \phi ) - (1+h-2\rho) ] r 
\le \alpha' h r  \text{ thanks to (\ref{C21})} \label{majax},
\end{eqnarray}
and, on the other hand,
\begin{eqnarray} 
\norme{z-v_z}{p_2} & \ge & \norme{z}{p_2}-\norme{v_z}{p_2} \ge 2 \rho r  \label{modulezvMin}.
\end{eqnarray}

\noindent
\underline{Geometrical fact:} 
Let us see that, for every $z \in \shell(T, s+(1+h)r, s+(1+h)^2 r)$,
\begin{equation}
\label{fact}
\mathcal{B}_{p_2}^{v_z}( \rho \norme{z-v_z}{p_2}) \subset \mathcal{B}_{p_2}^{z}( (1+\rho)  \norme{z-v_z}{p_2}) \cap \shell(S  \oplus \frac{\varphi r}{2(s+r)}, r, (1+h)r).
\end{equation}
The triangle ensures the first inclusion 
$\mathcal{B}_{p_2}^{v_z}( \rho \norme{z-v_z}{p_2}) \subset \mathcal{B}_{p_2}^{z}( (1+\rho)  \norme{z-v_z}{p_2}).
$
Let then $u \in \mathcal{B}_{p_2}^{v_z}( \rho \norme{z-v_z}{p_2} )$, then, by Lemma~\ref{direction},
\begin{eqnarray*} 
\norme{\frac{u}{\norme{u}{p_2}}-\chapeau v_z}{p_2} & \le & \frac{2\rho \norme{z-v_z}{p_2}}{\norme{v_z}{p_2}}  \\
& \le & \frac{2\rho \alpha' hr}{s+(1+h-2\rho)r} \text{ by Equation~(\ref{majax}) and
  definition of } v_z \\
& \le & \frac{\phi r}{2(s+r)} \text{ thanks to Equation~(\ref{C22})}
\end{eqnarray*}
and thus $\frac{u}{\norme{u}{p_2}}\in S \oplus \frac{\phi r}{2(s+r)}$.
For the norm of $u$, by definition of $v_z$ and Equation~(\ref{majax}), we have:
\begin{eqnarray*}
\norme{v_z}{p_2}- \rho \norme{z-v_z}{p_2} \le & \norme{u}{p_2} & \le  \norme{v_z}{p_2}+
\rho \norme{z-v_z}{p_2} \\ 
s+(1+h-2\rho)r - \rho \alpha' h r \le &  \norme{u}{p_2} & \le  s+(1+h-2\rho)r +
\rho \alpha'  h r  \\
s+r \le & \norme{u}{p_2} & \le s+(1+h)r, 
\end{eqnarray*}
thanks to Equations~(\ref{C24}) and (\ref{C21}).
This proves the geometrical fact~(\ref{fact}).

\vspace{0.2cm}
\noindent
\underline{Probabilistic estimate:} 
We can then estimate the probability of $E^c$.
Note that, with~(\ref{fact}), we have
\begin{eqnarray*}E^c \subset
\bigcup
&&
\left\{
\mathcal B_{p_2}^z((1+\rho)\norme{z-v_z}{p_2}) \not\subset B_{p_2}^z((1+\rho)(1+\epsilon)\norme{z-v_z}{p_2})
\right\} \\
& \cup & \{ B_{p_2}^z((1+\rho)(1+\epsilon)\norme{z-v_z}{p_2}) \not\subset \mathcal B_{p_2}^z((1+\rho)(1+\epsilon)^2\norme{z-v_z}{p_2}) \},
\end{eqnarray*}
where the union is for $z \in \shell(T, s+(1+h)r, s+(1+h)^2 r)$.
By Proposition~\ref{shapeGD}, there exist two strictly positive 
constants $A_2$ and $B_2$ such that 
for every $S$, for every $s,r>0$,
\begin{eqnarray*}
\P(E^c)
& \le & \sum_{ z \in \shell(T, s+(1+h)r, s+(1+h)^2 r) }
A_2 \exp (-B_2  (1+\rho) (1+\epsilon)\norme{z-v_z}{p_2})\\
& \le & \Card{\shell(T, s+(1+h)r, s+(1+h)^2 r)} \times A_2 \exp (-B_2 (1+\rho)2 \rho r) 
\end{eqnarray*}
thanks to (\ref{modulezvMin}). Then,
for every $z \in \shell(T, (1+h)r, (1+h)^2 r)$, thanks to (\ref{majax})
and (\ref{C25}), one has $(1+\epsilon)(1+\rho) 
\norme{z-v_z}{p_2}\le \alpha h r $, which proves the exponential estimate of the lemma. 

\vspace{0.2cm}
\noindent
\underline{Control of the infection paths:} 
It remains to estimate the minimal room needed to perform this infection, or in
other words to control
$$\bigcup_{z \in \shell(T, s+(1+h)r, s+(1+h)^2 r)} \mathcal{B}_{p_2}^{z}(
(1+\epsilon)^2(1+\rho) \norme{z-v_z}{p_2} ).$$ 
Let $z \in \shell(T , (1+h)r, (1+h)^2 r)$ and $u \in \mathcal{B}_{p_2}^{z}(
(1+\epsilon)^2(1+\rho) \norme{z-v_z}{p_2})$. We have:
\begin{eqnarray*}
\norme{u}{p_2} & \ge & \norme{z}{p_2} - (1+\epsilon)^2(1+\rho) \norme{z-v_z}{p_2}\\
& \ge & (1+\epsilon)^2(1+\rho)[s+(1+h-2\rho)r ] \\
&& -\left[(1+\epsilon)^2(1+\rho)\left(1+\frac{\phi (1+h)r}{s +(1+h)r} \right)-1\right]\norme{z}{p_2} 
\quad \quad \text{ thanks to (\ref{modulezvMaj})} \\
& \ge & (1+\epsilon)^2(1+\rho)[s+(1+h-2\rho)r ] \\
&&-
\left[(1+\epsilon)^2(1+\rho)\left(1+\frac{\phi (1+h)r}{s +(1+h)r} \right)-1\right](s+(1+h)^2r) \\
& \ge & s+ r [ (1+\epsilon)^2(1+\rho)   [ (1+h-2\rho)  -\phi -(1+\phi (1+h)r)(1+h)^2] +(1+h)^2     ].
\end{eqnarray*}
This last term tends to $s+r(1+h-\phi-(1+h)^3\phi)$
when $\epsilon$ and $\rho$ tend to $0$. By decreasing if
necessary $\epsilon$ and $\rho$, we obtain, as $h<1/2$:
$$\norme{u}{p_2} \ge s+(1-3\phi)(1+h)r.$$
Finally, by applying Lemma~\ref{direction} and then Inequality~(\ref{majax}), we have
$$
\norme{\chapeau u -\chapeau z}{p_2}
\le \frac{2\norme{u-z}{p_2}}{\norme{z}{p_2}} 
\le  \frac{2 (1+\epsilon)^2(1+\rho)\alpha' hr}{s+(1+h)^2r}
 \le \frac{2\alpha hr}{s+r}.
$$
Thus $u \in \shell(T\oplus\frac{2\alpha hr}{s+r}, s+(1-3\phi)(1+h)r,\infty)$, which ends the proof of the lemma.
\end{proof}

\begin{lemme}
\label{avantage2}
Let $R>0$ and $\eta \in (0,1/2)$. There exist two strictly positive constants $A,B$ such that 
$$
\forall\chapeau v \in \mathcal S \quad \forall t>0 \quad \P \left( \mathcal{G}^1\cap \Shadow(\chapeau v,t, Rt^{1/2+\eta}) \right)\le A\exp(-Bt^\eta).
$$
\end{lemme}

\begin{proof}
For convenience, we note $\hat x=v/\norme{v}{p_2}\in \mathcal S_{p_2}$. Therefore, $\chapeau{v}=\frac{\hat x}{\|\hat x\|_2}\in \mathcal S$.
Let $\delta>0$. 

\vspace{0.2cm}
\noindent
\underline{Idea of the proof:} 
The idea is quite natural: start the progression by the initialization Lemma~\ref{petit-bout-de-coquille}, and apply recursively the progression Lemma~\ref{progressionter} until the stronger infection surrounds the weaker one. The point is to ensure that this progression is not disturbed by the spread of the weaker infection.

\vspace{0.2cm}
\noindent
\underline{Step 0. Choice of constants:}
Choose $R>0$, $\eta \in (0,1/2)$ and $\delta \in (0,R/K_1)$. Choose then $\delta'>\delta$ such that $\delta'<\min \{ \delta/C_{p_1,p_2}, R/K_1 \}$. Lemma~\ref{petit-bout-de-coquille} gives then $\theta_0>0$ and $\gamma'_0>0$.
Remember that $C_{p_1,p_2}<1$ and choose, in this order, $\alpha,\gamma,\theta,h,\phi, \epsilon$:
\begin{align}
1&<\alpha  <\min \left\{\frac{1}{C_{p_1,p_2}},2\right\}, \label{3alpha} \\
1& <\gamma<\gamma'_0 \text{ and }0 < \theta<\theta_0 \label{m3gamma},\\
0& < h < \frac12 \text{ and } 1+h<\min\left\{\alpha,\frac{\gamma'_0}{\gamma}\right\} \label{m3h} \\
& \text{ and } 8\alpha h C_{p_1,p_2}C_{p_2,p_1}  <1  -\max \left\{\frac1\gamma,C_{p_1,p_2} \alpha \right\} \nonumber , \\
\phi&>0 \text{ and } (1+h)\left(1+\frac{3\phi}h \right)<\alpha \text{ and } \phi \gamma \delta' <2\theta\label{m3phi}\\
&\text{ and } \phi\left(4C_{p_1,p_2}C_{p_2,p_1} \left(\frac1{2h}+1\right)+3\right)+8\alpha h C_{p_1,p_2}C_{p_2,p_1}  \nonumber \\
&< 1  -\max \left\{\frac1\gamma,C_{p_1,p_2} \alpha \right\}, \nonumber\\
\epsilon&>0 \text{ and } \phi\left(4C_{p_1,p_2}C_{p_2,p_1}(\frac1{2h}+1)+3\right)+8\alpha h C_{p_1,p_2}C_{p_2,p_1}\label{3epsilon}\\ 
&  +\epsilon\alpha C_{p_1,p_2}<1 -\max \left\{\frac1\gamma,C_{p_1,p_2} \alpha \right\}\nonumber
. 
\end{align}
Set $\gamma'=(1+h)\gamma<\gamma'_0$.

\vspace{0.2cm}
\noindent
\underline{Step 1. Initialization of the spread:}
Let us introduce the following notations: 
\begin{align*}
E_1^1(\hat x, t) & =  \{ \eta^1(t+\delta t^{1/2+\eta})\, \subset \, \mathcal B_{p_1}^0( t+\delta't^{1/2+\eta}) \}, \\
E_1^2(\hat x, t) & =  \left\{  \eta^2(t+\delta t^{1/2+\eta}) \, \supset \, 
 \shell \left(\hat x\oplus
\frac{\phi r_{0}\norme{\hat x}{p_1}}{2(t+r_{0}\norme{\hat x}{p_1})}, \frac{t+\gamma \delta't^{1/2+\eta}}{\norme{\hat x}{p_1}}, \frac{t+\gamma' \delta't^{1/2+\eta}}{\norme{\hat x}{p_1}} \right) \right\}, \\  
E_1(\hat x, t) & =  E_1^1(\hat x, t) \cap E_1^2(\hat x, t).
\end{align*}
From Assumption $\phi\gamma\delta'<2\theta$ in~(\ref{m3phi}), it follows that $\frac{\phi r_{0}\norme{\hat x}{p_1}}{2(t+r_{0}\norme{\hat x}{p_1})}
\le \theta t^{-1/2+\eta}$ for $t$ large enough.
Then, by Lemma~\ref{petit-bout-de-coquille}, there exist two strictly positive
constants $A_1$ and $B_1$ such that for every $\hat x \in \mathcal{S}_{p_2}$, for every $t>0$, we have
\begin{equation}
\P\left(  \Shadow\left(\frac{\hat x}{\|\hat x\|_2}, t, Rt^{1/2+\eta}\right)
\backslash  E_1(\hat x, t) \right) \le  
A_1\exp(-B_1t^\eta). \label{D11}
\end{equation}
Thus, if $\Shadow(\frac{\hat x}{\|\hat x\|_2}, t, Rt^{1/2+\eta})$ occurs, then at the slightly larger
time $t_1(\hat x,t)=t+\delta t^{1/2+\eta}$, the first shell 
$$S_1(\hat x,t)= \shell \left(\hat x\oplus \frac{\phi r_{0}\norme{\hat x}{p_1}}{2(t+r_{0}\norme{\hat x}{p_1})} , \frac{t+\gamma \delta't^{1/2+\eta}}{\norme{\hat x}{p_1}}, \frac{t+\gamma' \delta't^{1/2+\eta}}{\norme{\hat x}{p_1}} \right)$$
is with high probability colonized by the $p_2$-infection. 
We want now to extend this colonization to larger and larger shells by
applying  recursively Lemma~\ref{progressionter}.

\vspace{0.2cm}
\noindent
\underline{Notations:}
We still need to introduce a certain number of notations, inspired by Lemma~\ref{progressionter}.
$$
\begin{array}{|l|}
\hline
k=1 \\
\hline
\displaystyle r_0=r_0(\hat x,t)= \frac{1}{\norme{\hat x}{p_1}}\gamma \delta't^{1/2+\eta}\\
r_1=r_1(\hat x,t)=(1+h)r_0\\
t_1=t_1(\hat x,t)= t+\delta t^{1/2+\eta}\\
\displaystyle A_1=A_1(\hat x,t)=\{\hat x\} \oplus \frac{\phi r_{0}\norme{\hat x}{p_1}}{2(t+r_{0}\norme{\hat x}{p_1})}  \\
\displaystyle S_1=S_1(\hat x,t)= \shell\left(A_1,\frac{t}{\norme{\hat x}{p_1}}+r_0,\frac{t}{\norme{\hat x}{p_1}}+r_1\right)\\
\hline
k \ge 2 \\
\hline
r_k=r_k(\hat x,t)=(1+h)r_{k-1}=(1+h)^{k}r_0
\quad \text{and} \quad r_k^{\min}=[1-3\phi](1+h)r_{k-2} \\
t_k=t_k(\hat x,t)=t_{k-1}+h \alpha  r_{k-2} = t+\delta t^{1/2+\eta}+
r_0\alpha[(1+h)^{k-1}-1] \\
\displaystyle A_k=A_k(\hat x,t)=A_{k-1}\oplus \frac{\phi r_{k-1}\norme{\hat x}{p_1}}{2(t+r_{k-1}\norme{\hat x}{p_1})} \\
\displaystyle A_k^+=A_k^+(\hat x,t)=A_{k-2}\oplus \frac{(\phi+2\alpha h) r_{k-2}\norme{\hat x}{p_1}}{t+r_{k-2}\norme{\hat x}{p_1}} \\
\displaystyle S_k=S_k(\hat x,t)=\shell \left(A_k,    \frac{t}{\norme{\hat x}{p_1}}+r_{k-1},\frac{t}{\norme{\hat x}{p_1}} +r_{k}\right) \\
\displaystyle S_k^+=S_k^+(\hat x,t)=\shell \left( A_k^+,  \frac{t}{\norme{\hat x}{p_1}}+ r_k^{\min} ,\infty \right) \\
\hline
\end{array}
$$
Define also the following events, for $k \ge 2$ and $x \in \Zd \backslash
\{0\}$:
\begin{align*}
E_k^1&=E_k^1 (\hat x,t) = \{\eta^1(t_k(\hat x,t)) \cap S_k^+(\hat x,t) =\varnothing \}, \\ 
E_k^2&=E_k^2 (\hat x,t) =  \{\eta^2(t_k(\hat x,t)) \supset S_k(\hat x,t) 
\}, \\ 
E_k&=E_k (\hat x,t) =  E_k^1 (\hat x,t) \cap E_k^2(\hat x,t).
\end{align*}
The aim is the following: we want to apply Lemma~\ref{progressionter} to prove that if $E_k^2 (\hat x,t)$ is
realized, then with high probability $E_{k+1}^2 (\hat x,t)$ is also realized. But we
need first to control the spread of the slow $p_1$-infection, and to see
that it will not disturb the spread of the fast $p_2$-infection from
$S_k(\hat x,t)$ to $S_{k+1}(\hat x,t)$.

\vspace{0.2cm}
\noindent
\underline{Step 2. Rough control of the slow $p_1$-infection:} 
Here, for convenience, the complementary event of $A$ is denoted by $\complement (A)$. Let $\epsilon>0$. Denote, for every $k \ge 2$
$$F_k^1 (\hat x,t) =\left\{\eta^1(t_k(\hat x,t)) \subset \mathcal B^0_{p_1}( 
t_1+(\delta'-\delta) t^{1/2+\eta}+(1+\epsilon)(t_k-t_1))\right\}.
$$
Let us prove that
there
exist two strictly positive constants $A_2$ and $B_2$ such that 
\begin{equation}
\label{EE2}
\forall t>0\quad\forall\hat x \in \mathcal S_{p_2}  \quad 
\P \left( E_1^1(\hat x,t)  \cap \complement \left( \miniop{}{\bigcap}{k \ge 2} F_k^1(\hat x,t) \right) \right) \le A_2 \exp (-B_2 t^{1/2+\eta}).
\end{equation} 
Note that
$$
\{B_{p_1}^0(s)\subset\mathcal{B}_{p_1}^0(s')\}\cap
\{B_{p_1}^0(s+t)\not\subset\mathcal{B}_{p_1}^0(s'+t')\}\subset
\miniop{}{\cup}{x\in\mathcal{B}^0_{p_1}(s')}\{ B_{p_1}^x(t)\not\subset\mathcal{B}_{p_1}^x(t')\},
$$
and thus that
$$\P \left( \{B_{p_1}^0(s)\subset\mathcal{B}_{p_1}^0(s')\}\cap
\{B_{p_1}^0(s+t)\not\subset\mathcal{B}_{p_1}^0(s'+t')\} \right)
\le \Card{\mathcal{B}_{p_1}^0(s')} \P \left(B_{p_1}^0(t)\not\subset\mathcal{B}_{p_1}^0(t') \right).
$$
In our context, this gives
\begin{eqnarray*}
&& \P \left( E_1^1(\hat x,t)  \cap \complement \left( \bigcap_{k \ge 2} F_k^1(\hat x,t) \right) \right) \\
& \le & \sum_{k \ge 2} \P \left( E_1^1(\hat x,t)  \cap \{ B_{p_1}^0(t_1+(t_k-t_1)) \not\subset
  \mathcal{B}_{p_1}^0(t_1+(\delta'-\delta) t^{1/2+\eta}+(1+\epsilon)(t_k-t_1) \} \right)  \\
& \le & \sum_{k \ge 2} \Card{\mathcal{B}_{p_1}^0(t+\delta' t^{1/2+\eta})} \P \left( {B}_{p_1}^0(t_k-t_1) \not\subset
\mathcal{B}_{p_1}((1+\epsilon)(t_k-t_1)) \} \right).
\end{eqnarray*}
The large deviation result, Proposition~\ref{shapeGD}, gives then two strictly positive constants $A,B$ such that
\begin{eqnarray*} \P \left( E_1^1(\hat x,t)  \cap \complement \left( \bigcap_{k \ge 2} F_k^1(\hat x,t) \right) \right) 
& \le & \Card{\mathcal{B}_{p_1}^0(t+\delta' t^{1/2+\eta})} \sum_{k \ge 2} A \exp \left( -B(t_k-t_1)  \right) \\ 
& \le &  \Card{\mathcal{B}_{p_1}^0(t+\delta' t^{1/2+\eta})}
A \sum_{k \ge 2} \exp \left( -B\alpha r_0 (k-1)h \right) \\
& \le & A \Card{\mathcal{B}_{p_1}^0(t+\delta' t^{1/2+\eta})} \frac{\exp \left( -B\alpha r_0 h \right) }{1-\exp \left( -B\alpha r_0 h \right)} \\
 & \le & A_2 \exp (-B_2 t^{1/2+\eta}),
\end{eqnarray*}
since $r_0=\frac1{\norme{\hat x}{p_1}}\gamma \delta't^{1/2+\eta}$.

\vspace{0.2cm}
\noindent
\underline{Step 3. Estimates for angles:}
Let us see that 
for any $\hat x \in \mathcal{S}_{p_2}$, for any $\phi, \psi\ge 0$,
$$(\hat x \oplus \phi) \oplus \psi \subset\hat x \oplus (\phi + \psi).$$
Let $z \in (\hat x \oplus \phi) \oplus \psi$: there exist $y \in\hat x \oplus \phi$ and $v\in \mathcal B_{p_2}(\psi)$ such that $z=y+v$.
As $y \in\hat x \oplus \phi$, there exists $w \in \mathcal B_{p_2}(\phi)$ such that $y=\hat x+w$. Thus
$$\norme{z-\hat x}{p_2} \le \norme{v}{p_2} +\norme{w}{p_2} \le \phi+\psi.$$
This implies
\begin{eqnarray*}
\forall k \ge 2 \quad A_k(\hat x,t)  & \subset &\hat x \oplus \left(\sum_{j=0}^{k-1} \frac{\phi r_j\norme{\hat x}{p_1}}{2(t+r_j\norme{\hat x}{p_1})} \right)
\text{ and } A_k^+(\hat x,t) \subset \hat x \oplus\theta_k^+, \\
 & &\text{with }\theta_k^+=\sum_{j=0}^{k-2} \frac{\phi r_j\norme{\hat x}{p_1}}{2(t+r_j\norme{\hat x}{p_1})} +\frac{(\phi+2\alpha h)  r_{k-2}\norme{\hat x}{p_1}}{t+r_{k-2}\norme{\hat x}{p_1}}.
\end{eqnarray*}
Then, for all $k \ge 2$, for $t$ large enough,
\begin{eqnarray}
\frac{t}{\norme{\hat x}{p_1}}\theta_k^+ & = & \frac{t}{\norme{\hat x}{p_1}} \left(\sum_{j=0}^{k-2} \frac{\phi r_j\norme{\hat x}{p_1}}{2(t+r_j\norme{\hat x}{p_1})} +\frac{(\phi+2\alpha h)  r_{k-2}\norme{\hat x}{p_1}}{t+r_{k-2}\norme{\hat x}{p_1}} \right) \nonumber \\
& \le & \frac{\phi}{2}\sum_{j=0}^{k-2}r_j\norme{\hat x}{p_1}+(\phi+2\alpha h) r_{k-2} \nonumber \\
& \le & r_0\left( (1+h)^{k-1}\left(\frac\phi{2h}+\phi+2\alpha h\right)\right)
. \label{majangle}
\end{eqnarray}

\vspace{0.2cm}
\noindent
\underline{Step 4. The weak can not bother the strong:} \\
Let us see that for all $k \ge 2$
\begin{equation}
\label{3bother}
\mathcal B^0_{p_1}( t+\delta' t^{1/2+\eta}+(1+\epsilon)r_0\alpha[(1+h)^{k-1}-1]) \cap S_k^+(\hat x, t)=\varnothing. 
 \end{equation}
Let $y$ such that $ \norme{y}{p_1} = t+\delta' t^{1/2+\eta}+(1+\epsilon)r_0\alpha[(1+h)^{k-1}-1])$ and $\norme{\hat y -\hat x }{p_2} \le \theta_k^+$. As $y=\norme{y}{p_1}\hat y/\norme{\hat y}{p_1}$, for $t$ large enough
\begin{eqnarray*}
\norme{y}{p_2} & \le &
 \frac{ t+\delta' t^{1/2+\eta}}{\norme{\hat x}{p_1}} 
+(t+\delta' t^{1/2+\eta}) \left|  \frac{1}{\norme{\hat y }{p_1}}-\frac1{\norme{\hat x}{p_1}}\right|
+\frac{(1+\epsilon)\alpha r_0 [(1+h)^{k-1}-1]}{\norme{\hat y }{p_1}} \\
& \le &
\frac{ t+\delta' t^{1/2+\eta}}{\norme{\hat x}{p_1}} 
+\frac{(t+\delta't^{1/2+\eta}) C_{p_2,p_1} \norme{\hat y -\hat x }{p_2} }{\norme{\hat x}{p_1}\norme{\hat y }{p_1}}
+\frac{(1+\epsilon)\alpha  r_0[(1+h)^{k-1}-1]}{\norme{\hat y }{p_1}} \\
& \le &
\frac{ t}{\norme{\hat x}{p_1}}  +  \frac{r_0}{\gamma} 
+\frac{2t C_{p_1,p_2}C_{p_2,p_1} \theta_k^+ }{\norme{\hat x}{p_1}}
+C_{p_1,p_2}{(1+\epsilon)\alpha r_0[(1+h)^{k-1}-1]} \\
& \le &
\frac{ t}{\norme{\hat x}{p_1}}  +  \frac{r_0}{\gamma} 
+2 C_{p_1,p_2}C_{p_2,p_1}r_0 (1+h)^{k-1}\left(\frac\phi{2h}+\phi+2\alpha h\right)\\
& & +C_{p_1,p_2}{(1+\epsilon)\alpha r_0[(1+h)^{k-1}-1]}, 
\end{eqnarray*}
where the last inequality follows from~(\ref{majangle}). Then
\begin{eqnarray}
&& \frac{1}{r_0} \left(\norme{y}{p_2}-\frac{t}{\norme{x}{p_1}}-r_k^{\min} \right) \nonumber \\
& \le & \frac{1}{\gamma} 
+2 C_{p_1,p_2}C_{p_2,p_1}(1+h)^{k-1}\left(\frac\phi{2h}+\phi+2\alpha h\right) \nonumber \\
& &+C_{p_1,p_2}{(1+\epsilon)\alpha [(1+h)^{k-1}-1]}-(1-3\phi)(1+h)^{k-1} \nonumber\\
& \le & (1+h)^{k-1} \left(2 C_{p_1,p_2}C_{p_2,p_1}\left(\frac\phi{2h}+\phi+2\alpha h\right)+C_{p_1,p_2}(1+\epsilon)\alpha -1+3\phi \right) \nonumber\\
&& +\left( \frac{1}{\gamma} -C_{p_1,p_2}(1+\epsilon)\alpha\right) \label{3berk} .
\end{eqnarray}
We want to prove that this quantity is negative for every $k \ge 2$.
Conditions~(\ref{m3phi}) and (\ref{3epsilon}) ensures that the coefficient in $(1+h)^{k-1}$ in (\ref{3berk}) is negative. 
Thus, asymptotically in $k$, the right-hand side in (\ref{3berk}) is negative. To ensure it is negative for every $k\ge 2$, we only need to see that it is true for $k=2$, which is ensured by  Conditions~(\ref{m3phi}) and (\ref{3epsilon}).
This proves~(\ref{3bother}). Note that this also implies
\begin{equation}
\label{3EF}
\forall k \ge 2 \quad F_k^1 \subset E_k^1.
\end{equation}
Equations~(\ref{3EF}) and (\ref{EE2}) together give:
\begin{equation}
\label{3EE2}
\forall t>0 \quad \forall\hat x \in \mathcal{S}_{p_2}  \quad 
\P \left( E_1^1(\hat x,t)  \cap \left( \miniop{}{\bigcup}{k \ge 2} E_k^1(\hat x,t)^c \right) \right) \le A_2 \exp (-B_2 t^{1/2+\eta}).
\end{equation} 

\vspace{0.2cm}
\noindent
\underline{Step 5. Control of the fast $p_2$-infection:}
Let $A_3$ and $B_3$ be the two strictly positive constants given by Lemma~\ref{progressionter} for the choice for $\alpha, h, \phi$ we made in~(\ref{3alpha}), (\ref{m3h}) and (\ref{m3phi}). Note that for every $k\ge 2$, we have
$$E_k^1\cap E_{k-1}^2\cap E(A_{k-1},r_{k-2},\frac{t}{\norme{\hat x}{p_1}})\subset E_k^2$$
where the event $E(.,.,.)$ was defined in Lemma~\ref{progressionter}.
Thus, the application of  Lemma~\ref{progressionter} implies that
for 
any $\hat x \in \mathcal{S}_{p_2}$, any $t>0$, for every $k\ge2$,
\begin{eqnarray*}
\P((E_k^2)^c \cap E_k^1\cap E_{k-1}^2)
& \le & \P(E(A_{k-1},r_{k-2},\frac{t}{\norme{\hat x}{p_1}})^c) \le A_3 (\frac{t}{\norme{\hat x}{p_1}}+r_{k-2})^d\exp(-B_3r_{k-2}).
\nonumber
\end{eqnarray*}
Thus, for each $t\ge 1$, each $\hat x\in\mathcal{S}_{p_2}$,
\begin{eqnarray}
\sum_{k \ge 2} \P \left( (E_k^2)^c\cap (E_k^1\cap
  E_{k-1}^2\right) 
& \le & A_3 \sum_{k \ge 2} \left( \frac{t}{\norme{\hat x}{p_1}}+r_{k-2}\right)^d\exp(-B_3r_{k-2} ) \nonumber \\
& \le & A_3 \sum_{k \ge 0} ( C_{p_1,p_2}t+1)^d(r_{k}+1)^d\exp(-B_3r_{k} ) \nonumber \\
& \le & A_4 ( C_{p_1,p_2}t+1)^d\sum_{k \ge 0} \exp(-B_4r_{k}) \nonumber \\
& \le & A_4 ( C_{p_1,p_2}t+1)^d\sum_{k \ge 0} \exp(-B_4r_{0}(1+kh)) \nonumber \\& \le & A_4 ( C_{p_1,p_2}t+1)^d\exp(-B_4r_{0})(1- \exp(-B_4r_{0}h))^{-1}\nonumber\\
&  \le & A_5 \exp(-B_5 t^{1/2+\eta}).\label{EE5} 
\end{eqnarray}
where $A_4,A_5$ and $B_4,B_5$ are strictly positive constants.

\vspace{0.2cm}
\underline{Conclusion:} 
For $k$ large enough, the set $S_k$ disconnects
$0$ from 
infinity, and thus the event $\bigcap_{k\ge1} E_k$ implies that the slow
$p_1$-infection is surrounded by the fast $p_2$-infection and thus
dies out. So, using (\ref{D11}), (\ref{3EE2}) and (\ref{EE5}), we obtain:
\begin{eqnarray*}
&& \P\left(\mathcal{G}^1\cap \Shadow\left(\frac{\hat x}{\|\hat x\|_2}, t,Rt^{1/2+\eta}\right)\right) \\
& \le & \P \left( \Shadow\left(\frac{\hat x}{\|\hat x\|_2}, t,Rt^{1/2+\eta}\right)\cap \bigcup_{k\ge 1} E_k(\hat x, t)^c \right) \\
& \le & \P \left( \Shadow\left(\frac{\hat x}{\|\hat x\|_2}, t,Rt^{1/2+\eta}\right) \cap E_1(\hat x, t)^c \right)
+ \P \left( E_1^1(\hat x,t)  \cap \bigcup_{k\ge 2} E_k^1(\hat x, t)^c \right) \\
&& + \sum_{k \ge 2} \P \left( (E_k^2(\hat x, t))^c\cap (E_k^1(\hat x, t) \cap
  E_{k-1}^2(\hat x, t)\right) \\
& \le & A\exp(-Bt^\eta),
\end{eqnarray*}
which completes the proof.
\end{proof}

\begin{proof}[Proof of Theorem~\ref{avantage3}]
Proposition~\ref{shapeGD} and Lemma~\ref{coupling} give the existence
of strictly positive constants $\alpha,\beta,A_1,B_1$ such that the event
$$F_t=\{\forall y\in \partial\eta(t)\quad \|y\|_2\in (\alpha t,\beta t)\}$$
has a probability larger than $1-A_1\exp(-aB_1t)$,
so we only have to control the probability of the event
$$\mathcal{G}^1\cap \Shade(t,Mt^{1/2+\eta})\cap F_t.$$
Assume that $t$ is large enough to have $\alpha t>Mt^{1/2+\eta}$ and set
$$\theta=\min\left\{\frac12,\left(\frac{Mt^{1/2+\eta}}{2\beta t}\right)^2\right\}.$$
By Lemma~\ref{familledeboules},  there exists a subset $T$ of the unit sphere
$\mathcal S$ with $\Card T \le C(1+\frac1{\theta})^{d-1}$ such that
$\mathcal S\subset \miniop{}\cup{\chapeau x\in T}\mathcal{B}_2^{\chapeau x}(\theta)$.

Assume now that $ \Shade(t,Mt^{1/2+\eta}) \cap F_t$ occurs: there exists $\chapeau u\in\mathcal S_2$ such that  
$\Shadow(\chapeau u,t,Mt^{1/2+\eta})$ happens. Let $\chapeau x\in S$ be such that $\|\chapeau u-\chapeau x\|_2\le\theta$. Let us prove that 
$\Shadow(\chapeau x,t, \frac{M}2t^{1/2+\eta})$ happens.

Let $\gamma$ be an infinite path in $\Cyl_+(\chapeau x,\frac{M}2t^{1/2+\eta})$ starting at some point $y\in\eta(t)$. 
We must prove that $\gamma$ meets $\partial\eta(t)\cap \eta^2(t)$.
We can suppose without loss of generality that  $y$ is the last point
of $\gamma$ in $\eta(t)$ and thus $y\in\partial\eta(t)$.
Note that the points in $\gamma$ after $y$ are in $\mathcal{B}_2(\alpha t)^c \subset \mathcal{B}_2(Mt^{1/2+\eta})^c$.
\begin{itemize}
\item Either $y\in\eta^2(t)$, and we are done.
\item Or $y\in\eta^1(t)$. Let $z$ be the point after $y$ along $\gamma$ where $\gamma$ exits from 
$\mathcal{B}_2(\beta t)$: between $y$ and $z$, the path  is in
$\Cyl_{+}(\chapeau x,\frac{M}2t^{1/2+\eta})\cap \mathcal{B}_2(\beta t)\cap\mathcal{B}_2(Mt^{1/2+\eta})^c\subset 
\Cyl_+(\chapeau u,Mt^{1/2+\eta})\cap\mathcal{B}_2(\beta t)$ thanks to Lemma~\ref{petitcyl}. We build now a path $\gamma'$ inside 
the $ \Cyl_+(\chapeau u,Mt^{1/2+\eta})$ by concatenating the portion of $\gamma$ between $y$ and $z$ and any infinite path starting from $z$ and staying in  $ \Cyl_+(\chapeau u,Mt^{1/2+\eta}) \cap \mathcal{B}_2(\beta t)^c$: this path prevents the occurrence of $\Shadow(\chapeau u,t,Mt^{1/2+\eta})$, as it starts from a point in $\eta^1(t)$ and do not visit any other point in $\eta(t)$. The assumption $y\in\eta^1(t)$ is thus contradicted.
\end{itemize}
So $y\in\eta^2(t)$, which means that $\Shadow(\chapeau x,t, \frac{M}2t^{1/2+\eta})$ happens, and implies
$$ \Shade(t,Mt^{1/2+\eta}) \cap F_t \subset\miniop{}\cup{\chapeau x\in T}  \Shadow\left(\chapeau x,t,\frac{M}2t^{1/2+\eta}\right).$$
Finally, Theorem~\ref{avantage2} give two strictly positive constants $A',B'$ such that
\begin{eqnarray*}  
\P( \mathcal{G}^1\cap \Shade(t,Mt^{1/2+\eta})\cap F_t)
&\le&  \Card T  A'\exp(-B't^{\eta})\\ 
& \le & A'\left(1+\frac1{\theta}\right)^{d-1}\exp(-Bt^{\eta})\\
& \le & A'\left(1+\frac{4 \beta^2}{M^2} t^{1-2\eta}\right)^{d-1}\exp(-Bt^{\eta}),
\end{eqnarray*}
which ends the proof.
\end{proof}

\begin{proof}[Proof of Corollary~\ref{avantage4}]
We can assume that $\eta \in (0,1)$. For every $n \ge 1$, we define $t_n=n^2$ and
$$A(n,M)=\miniop{}{\cup}{t\in [t_n,t_{n+1}]} \Shade(t,Mt^{1/2+\eta}).$$
For each $M>0$, we have to prove that $\P \left(\mathcal{G}^1\cap \miniop{}{\limsup}{n\to +\infty}A(n,M) \right)=0$. By the Borel-Cantelli lemma, it is sufficient to 
prove that 
$$\sum_{n=1}^{+\infty}\P (\mathcal{G}^1\cap A(n,M))<+\infty.$$
Let then $M>0$. Obviously,
\begin{eqnarray*}
&& \P( \mathcal{G}^1\cap A(n,M)) \\
& \le & \P \left( \mathcal{G}^1\cap \Shade \left(t_{n+1},\frac M 2 t_{n+1}^{1/2+\eta} \right) \right)
+\P \left(\eta^1(t_{n+1})\not\subset\mathcal{B}_{p_1}(2t_{n+1}) \right)\\
& & +\P \left( \left(A(n,M)\backslash \Shade \left(t_{n+1},\frac M 2 t_{n+1}^{1/2+\eta}\right) \right)\cap\{\eta^1(t_{n+1})\subset\mathcal{B}_{p_1}(2t_{n+1})\} \right).
\end{eqnarray*}
Lemma~\ref{avantage3} gives two strictly positive constants $A_1,B_1$ such that $$\P \left( \mathcal{G}^1\cap \Shade \left(t_{n+1},\frac M 2 t_{n+1}^{1/2+\eta} \right) \right) \le A_1\exp(-B_1n^{2\eta}),$$ while
Proposition~\ref{shapeGD} gives two strictly positive constants $A_2,B_2$ such that  $$\P \left(\eta^1(t_{n+1})\not\subset\mathcal{B}_{p_1}(2t_{n+1}) \right) \le A_2\exp(-B_2n^2).$$ So it only remains to prove that
$$\sum_{n=1}^{+\infty}\P \left( \left(A(n,M)\backslash \Shade \left(t_{n+1},\frac M 2 t_{n+1}^{1/2+\eta}\right) \right)\cap\{\eta^1(t_{n+1})\subset\mathcal{B}_{p_1}(2t_{n+1})\} \right)<+\infty.$$
Assume now that $(A(n,M)\backslash \Shade(t_{n+1},\frac M 2 t_{n+1}^{1/2+\eta}))\cap\{\eta^1(t_{n+1})\subset\mathcal{B}_{p_1}(0,2t_{n+1})\}$ holds: there exists $\chapeau x$ and $t\in [t_n,t_{n+1})$ 
such that $\Shadow(\chapeau x,t,Mt^{1/2+\eta})$ holds but not $\Shadow(\chapeau x,t_{n+1},\frac M 2 t_{n+1}^{1/2+\eta})$.
Since $\Shadow(\chapeau x,t_{n+1},\frac M 2 t_{n+1}^{1/2+\eta})$ is not fulfilled, there exists some infinite path $\gamma$ in $\Cyl_+(\chapeau x, \frac M 2 t_{n+1}^{1/2+\eta})$  starting in some point $v\in\eta^1(t_{n+1})$ and that never meets $\eta^2(t_{n+1})$. By the definition of the process $(\eta^1(t))_{t\ge 0}$, there exist $u\in\eta^1(t_{n})$ and a path $\gamma'$ from $u$ to $v$
such that 
$$\sum_{e\in \gamma'}\omega_{e}^1\le t_{n+1}-t_n$$ 
and the path $\gamma'$ does not
meet any point in $\eta^2(\infty)$. The path $\gamma'$ can not stay completely in 
$\Cyl_+(\chapeau{x},Mt_n^{1/2+\eta})$, otherwise concatenating  $\gamma'$ and $\gamma$ together would contradict $\Shadow(\chapeau x,t,Mt^{1/2+\eta})$. Let $w$ be a point in $\gamma' \cap \Cyl_+(\chapeau{x},Mt_n^{1/2+\eta})^c$: the portion $\gamma''$ of $\gamma'$ between $w\in \eta^1(t_{n+1})\cap \Cyl_+(\chapeau{x},M t_n^{1/2+\eta})^c$ and $v\in \eta^1(t_{n+1})\cap \Cyl_+(\chapeau{x},\frac M 2 t_{n+1}^{1/2+\eta})$ satisfies
$$ \sum_{e\in \gamma''}\omega_{e}^1\le t_{n+1}-t_n.$$ 
Provided that $n$ is large enough, we can say
that $\norme{w-v}{p_1}\ge 2(t_{n+1}-t_n)$ -- note that $ t_n^{1/2+\eta} \gg t_{n+1}-t_n$.
Using the fact  that $\eta^1(t_{n+1})\subset\mathcal{B}_{p_1}(2t_{n+1})$, we obtain
\begin{eqnarray*}
 & & (A(n,M)\backslash \Shade(t_{n+1},M t_{n+1}^{1/2+\eta}/2))\cap\{\eta^1(t_{n+1})\subset\mathcal{B}_{p_1}(0,2t_{n+1})\}\\ 
 & \subset & \miniop{}{\cup}{w\in \mathcal{B}_{p_1}(2t_{n+1})}
\{B_{p_1}^w(t_{n+1}-t_n)\not\subset \mathcal{B}_{p_1}^w(2(t_{n+1}-t_n))\}.
\end{eqnarray*}
So it follows from Proposition~\ref{shapeGD} that there exist strictly positive constants $A_3,B_3$ such that
\begin{eqnarray*}
&&\P((A(n,M)\backslash \Shade(t_{n+1},M t_{n+1}^{1/2+\eta}/2))\cap\{\eta^1(t_{n+1})\subset\mathcal{B}_{p_1}(0,2t_{n+1})\}) \\
& \le&  \Card {\mathcal{B}_{p_1}(2t_{n+1})} A_3 \exp(-B_3(t_{n+1}-t_n))\\
& \le & K n^{2d}A_3 \exp(-2B_3n),
\end{eqnarray*}
which is summable.
This ends the proof of the theorem. 
\end{proof}

Note that the order $t^{1/2+\eta}$ that appears in our results follows from moderate deviations for fluctuations with respect to the asymptotic shape given in Proposition~\ref{shapeMD}. However, the conjectured order for these fluctuations is rather $1/3$. The proofs would apply with an estimate of the type
$$\forall t>0 \quad \P \left(
 \mathcal B^0 (t -\beta t^{1/3+\eta})\subset B^0(t) \subset  \mathcal B^0(t +\beta t^{1/3+\eta})
\right) \ge 1-  A \exp(-B t^\eta).$$
which would lead to replace $t^{1/2+\eta}$ by $t^{1/3+\eta}$ in our results.

\section{Moderate deviations for the global growth of the epidemics} 
\label{sectionmoderate}

This section aims to prove Theorem~\ref{pasdexcroissance}: when the weak survives, we recover the same fluctuations with respect to the asymptotic shape as in the case where the weak infection evolves alone.
 
In fact, the only point is to see that if the strong infection at time $t$ admits points outside
$\mathcal B_{p_1} \left( t + \beta t^{1/2+\eta} \right)$, then this positional advantage enables  to create an event of type $\Shade(T,MT^{1/2+\eta})$ 
at the slightly larger time $T=t+\alpha t^{1/2+\eta}$, and the probability of such an event is controlled by Lemma~\ref{avantage3}.

\begin{proof} 
Let $\beta>0$ and $\eta \in (0, 1/2)$. We want to prove that there exist two strictly positive constants $A,B$ such that
$$
\forall t >0 \quad \P \left(\mathcal G^1 \cap \left\{ \mathcal B_{p_1} \left( t - \beta t^{1/2+\eta} \right) \subset \eta(t) \subset \mathcal B_{p_1} \left( t + \beta t^{1/2+\eta} \right) \right\}^c \right) \le A \exp(-Bt^\eta).
$$ 
We can first easily rule out the cases where $ \eta^1(t)\not\subset \mathcal B_{p_1} \left( t + \beta t^{1/2+\eta} \right)$ and where $\mathcal B_{p_1} \left( t - \beta t^{1/2+\eta} \right) \not\subset \eta(t)$, by Proposition~\ref{shapeMD} and Lemma~\ref{coupling}: there exist two strictly positive constants $A_1,B_1$ such that for every $t >0$, we have 
\begin{eqnarray*}
& & \P \left(\mathcal G^1 \cap \left\{ \mathcal B_{p_1} \left( t - \beta t^{1/2+\eta} \right) \subset \eta(t) \subset \mathcal B_{p_1} \left( t + \beta t^{1/2+\eta} \right) \right\}^c \right) \\
& \le & 
\P \left( \mathcal B_{p_1} \left( t - \beta t^{1/2+\eta} \right) \not\subset \eta(t)\right)+ \P \left( \eta^1(t)\not\subset \mathcal B_{p_1} \left( t + \beta t^{1/2+\eta} \right) \right)\\
& &+  \P \left(\mathcal G^1 \cap \left\{\eta^2(t)\not\subset \mathcal B_{p_1} \left( t + \beta t^{1/2+\eta} \right) \right\}\right)   \\
& \le & \P \left(\mathcal G^1 \cap \left\{\eta^2(t)\not\subset \mathcal B_{p_1} \left( t + \beta t^{1/2+\eta} \right) \right\}\right)+ A_1\exp(-B_1t^{\eta}).
\end{eqnarray*}
Remember that $C_{p_2,p_1}$ is a strictly positive constant such that
$$\forall x\in\Rd\quad \norme{x}{p_1} \le C_{p_2,p_1} \norme{x}{p_2}.$$
Choose now $\alpha>0$ such that  $(2C_{p_2,p_1}+3)\alpha<\beta$. 
Define now the following events:
\begin{eqnarray*}
E & = & \left\{ \forall x\in \mathcal B_{p_1}(2t)\quad   \mathcal{B}^x_{p_2} \left(\frac{\alpha}2 t^{1/2+\eta}\right) \subset
{B}^x_{p_2}\left(\alpha t^{1/2+\eta}\right) \subset
\mathcal{B}^x_{p_2} \left(2\alpha t^{1/2+\eta}\right)
\right\}, \\
F & = & \left\{ \eta^2(t) \not\subset \mathcal B_{p_1} \left( t + \beta t^{1/2+\eta} \right) \right\} \cap 
\left\{ \eta^1(t+\alpha t^{1/2+\eta}) \subset \mathcal B_{p_1} \left( t + 3\alpha t^{1/2+\eta} \right) \right\}.
\end{eqnarray*}
By Proposition~\ref{shapeGD} applied to the $p_2$-epidemic, there exist positive constants $A_2,A'_2,B_2,B'_2$ such that 
$$\P(E^c) \le \Card{\mathcal B_{p_1}(2t)} A'_2\exp(-B'_2t^{1/2+\eta})\le A_2 \exp(-B_2t^\eta),$$
and, by Proposition~\ref{shapeMD} applied to the $p_1$-epidemic, there exist positive constants $A_3,B_3$ such that 
\begin{eqnarray*}
\P \left(  \left\{\eta^1(t+\alpha t^{1/2+\eta})\not\subset \mathcal B_{p_1} \left( t + 3\alpha t^{1/2+\eta} \right) \right\}\right)
&   \le & A_3\exp(-B_3t^{\eta}).
\end{eqnarray*}
It is thus sufficient to prove that
$$\P(\mathcal G_1 \cap E \cap F) \le A\exp(-Bt^\eta).$$
Assume now that $E\cap F$ happens: there exists $x$
with $t+\beta t^{1/2+\eta}\le\norme{x}{p_1}\le 2t$ and $x\in\eta^2(t)$. 
Note first that $\mathcal{B}_{p_1}(t+3\alpha t^{1/2+\eta}) \cap \mathcal{B}^x_{p_2}(2\alpha t^{1/2+\eta})=\varnothing$: if $z \in \mathcal{B}^x_{p_2}(2\alpha t^{1/2+\eta})$, then, by the choice we made for $\alpha$,
\begin{eqnarray*}
\norme{z}{p_1} & \ge & \norme{x}{p_1} -  \norme{z-x}{p_1} \ge t+\beta t^{1/2+\eta} -C_{p_2,p_1}\norme{z-x}{p_2} \\
& \ge & t+\beta t^{1/2+\eta} -2C_{p_2,p_1}\alpha t^{1/2+\eta} > t+3\alpha t^{1/2+\eta}.
\end{eqnarray*}
Now, by the event $F$, we have
$ \eta^1(t+\alpha t^{1/2+\eta}) \subset \mathcal{B}_{p_1}(t+3\alpha t^{1/2+\eta})$
and, by the event $E$, we also have
$B^x_{p_2}(\alpha t^{1/2+\eta}) \subset
\mathcal{B}^x_{p_2}(2\alpha t^{1/2+\eta}).$
Thus, as $x \in\eta^2(t)$ and $\mathcal{B}_{p_1}(t+3\alpha t^{1/2+\eta}) \cap \mathcal{B}^x_{p_2}(2\alpha t^{1/2+\eta})=\varnothing$, we obtain that
$$\eta^2(t+\alpha t^{1/2+\eta})\supset  B^x_{p_2}(\alpha t^{1/2+\eta}) \supset \mathcal{B}^x_{p_2}\left(\frac\alpha2 t^{1/2+\eta}\right).$$
Let $C>0$ be such that $\forall x\in\Rd, \;\|x\|_2\ge C\norme{x}{p_2}$. Define $M=\frac{C\alpha}{2(1+\beta)^{1/2+\eta}}$ and $T=t+\alpha t^{1/2+\eta}$. Then, for every $t \ge 1$, 
$$\mathcal{B}^x_{p_2}\left(\frac\alpha2 t^{1/2+\eta}\right) \supset \mathcal B_2^x(M  T^{1/2+\eta}).$$
This, with the previous inclusion, implies that $\eta^2(T)\supset \mathcal B_2^x(M  T^{1/2+\eta})$, then that $\eta^2(T)$ disconnects $\eta^1(T)$ from infinity
in $\Cyl_+(x/\|x\|_2,MT^{1/2+\eta})$, then that $\Shade(T,MT^{1/2+\eta})$ occurs. Then
$$\P(\mathcal G_1 \cap E \cap F)  \le \P(\mathcal G_1 \cap \Shade(T,MT^{1/2+\eta}) \le A\exp(-Bt^\eta)$$
by Lemma~\ref{avantage3}.
\end{proof}

\section{Density of the strong in the two dimensional case}
\label{sectiondimdeux}

This section is devoted to the proof of Theorem~\ref{th-densite}: we prove that in dimension two, when coexistence occurs, the strong epidemic finally occupies a subset of $\Z^2$ with null density. We first need some definitions:

\begin{defi} For any $t \ge 0$, denote by $C_{\text{ext}} (t)$ the infinite connected component of $\eta(t)^c$, and by $\partial^\infty \eta(t)$ the external boundary of $\eta(t)$:
$$\froufrou \eta(t) = \{z \in \eta(t): \; \exists y \in C_{\text{ext}} (t), \; y \sim z \}.$$
\end{defi}

\begin{figure}[h!]
\label{hollywood}
\scalebox{0.5}{\input{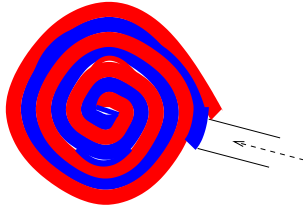}}
\caption{The shadow created by fast/blue on slow/red has a small diameter, while fast/blue fills a set with positive density.} 
\end{figure}

What is specific to the two dimensional case is that the external boundary of the fast infection $\froufrou \eta(t) \cap \eta^2(\infty)$ is $*$-connected -- this will be proved in Lemmas~\ref{connexdeux} and  \ref{connexeetoile}.  Loosely speaking, by Theorem~\ref{pasdexcroissance}, the external boundary $\froufrou \eta(t)$ of the infection at time $t$ is included in a very thin annulus with radius $t$ and width $t^{1/2+\eta}$. Then, Theorem~\ref{avantage4}, combined with the $*$-connectivity of $\froufrou \eta(t) \cap \eta^2(\infty)$, ensures that
the shadow cast by the fast infection  on the slow infection has a diameter smaller than $t^{1/2+\eta}$. However, it remains to control the points in $\eta^2(\infty)$ that are never in a position to create shadow, see Figure~\ref{hollywood}.

The proof breaks down in higher dimension, as we can imagine a configuration where the fast infection occupies a tree whose branches simultaneously widen and ramify. For instance, we can assume that the radius at height $t$ is of order $t^{1/2}$ and the number of branches at height $t$ is of order $t^{d-3/2}$.

Let us now recall the graphical duality of the square lattice.
Let $\Zp=\Z^2+(1/2,1/2)$, $\Edeuxstar=\{\{a,b\}: \; a,b\in\Zp\text{ and }\|a-b\|_2=1\}$ and $\Ldeuxstar=(\Zp,\Edeuxstar)$, which is isomorphic to $\Ldeux$.
For each bond $e=\{a,b\}$ of $\Ldeux$ (\resp $\Ldeuxstar$), let us denote by
$s(e)$ the only subset $\{i,j\}$ of $\Zp$ (\resp $\Z^2$) such that the quadrangle $aibj$ is 
a square in $\R^2$. The application $s$ is clearly an involution. 

For any finite set $A\subset\Zdeux$, we denote by $\Peierls(A)$ the set of Peierls contours associated 
to $A$, that is
$$\Peierls(A)=\{e\in\Edeuxstar: \;\1_A\text{ is not constant on }s(e)\}.$$
If, on the plane $\R^2$, we draw the edges which are in $\Peierls(A)$, we obtain a family of curves -- the so-called Peierls contours -- which are exactly
the boundary of the subset $A+[-1/2,1/2]^2$ of $\R^2$.
If $A\subset\Zdeux$ is a bounded connected subset of $\Ldeux$, there exists a unique set of bonds $\Gamma(A)\subset\Peierls(A)\subset \Edeuxstar$ which form a cycle  surrounding $A$, in the sense that
every infinite connected subset of bonds $D\subset\Edeux$ satisfying 
$D\cap A\ne\varnothing$ also satisfies $D\cap s(\Gamma(A))\ne\varnothing$.
If we draw  $\Gamma(A)$ on the plane $\R^2$, we get the external boundary
of $A+[-1/2,1/2]^2$, \ie the boundary of the infinite connected component
of $(A+[-1/2,1/2]^2)^c.$

Note also that if $\gamma$ is a Jordan curve on $\Ldeuxstar$, the set
$\text{Int}(\gamma)$ ($\text{Ext}(\gamma)$) composed by the points
in $s(\gamma)$ that are in the bounded (respectively,  unbounded) connected component
of $\R^2\backslash\gamma$ is $*$-connected.

We begin with two lemmas to prove the $*$-connectivity of the set $\eta^2(t) \cap \froufrou \eta(t)$
in dimension 2.

\begin{lemme}
\label{connexdeux}
Let  $A$, $B$ be two disjoint finite connected subsets of $\Z^2$ such that $A$, $B$, and $A\cup B$ are connected. We define
\begin{eqnarray*}
\Delta_A=\{e\in\Edeuxstar: \; s(e)\cap A\ne\varnothing\} & \text{ and } & 
\Delta_B=\{e\in\Edeuxstar: \; s(e)\cap B\ne\varnothing\},\\
E_{A \cup B}^A=\Gamma(A\cup B)\cap\Delta_A & \text{ and } & E_{A \cup B}^B=\Gamma(A\cup B)\cap\Delta_B.
\end{eqnarray*}
Then $E_{A \cup B}^A$ and $E_{A \cup B}^B$ are connected.
\end{lemme}

\begin{proof}
Since $\Gamma(A\cup B)=E_{A \cup B}^A \cup E_{A \cup B}^B$ is connected, we can assume without
loss of generality that $E_{A \cup B}^A$ and $E_{A \cup B}^B$ are non-empty.
The contour $\Gamma(A\cup B)$ is  a cycle that we denote as a sequence of distinct consecutive edges $e_0,e_1,\dots,e_{f-1}$, where
the ``end'' of  $e_{f-1}$ is the ``start'' of $e_0$.
To these edges we associate a sequence $x_0,\dots, x_{f-1}$ of points in $\Zdeux$
  such that $(A\cup B)\cap s(e_i)=\{x_i\}$.
We also denote by $m_0,\dots, m_{f-1}$ the middle points of the edges $e_0,e_1,\dots,e_{f-1}$.

Assume that $x_0\in A$ and suppose by contradiction that $E_{A \cup B}^A$ is not connected: there exist $p,q$ with 
$1<p<q< f$ with
 $x_i\in A$ for $i\in\{0\}\cup\{p,\dots,q-1\}$
whereas  $x_i\in B$ for $i\in \{1,\dots,p-1\}\cup\{q\}$.

Since $A$ is connected, there exists a simple path in $A$ from
$x_0$ to $x_{q-1}$ which corresponds to a path
in $\Zdeuxstar$ from the start of $e_0$ to the end of $e_{q-1}$. 
The union of this path with the path $(e_0,e_1,\dots e_{q-1})$
makes a Jordan curve $\gamma$. 

Obviously, $m_q\notin\gamma$. Since $m_q$ is on the outer boundary of the connected set $(A\cup B)+[-1/2,1/2]^2$, there exists
an (infinite) path in $(A\cup B)^c$ joining $m_q$ to infinity.
So, we can say that $m_q$ is in the infinite component of  $\gamma^c$.
Since $B$ is connected, there exists a path $\gamma'$ in $B$ from $x_q$ to $x_1$.
Let $e$ be the  first edge of  $\gamma'$ which crosses $\gamma$. 

By construction, we know that each edge $e\in\Ldeux$ which crosses $\gamma'$
from the unbounded component to a point in $B$ must be one of the $s(e_i)$'s.
But there is a contradiction because no $s(e_i)$ can have
both  ends in $A\cup B$.
\end{proof}

\begin{lemme}
\label{connexeetoile}
In dimension 2, the set $\eta^2(t) \cap \froufrou \eta(t) $ is $*$-connected. 
\end{lemme}

\begin{proof}
By the very definition of the evolution process, $\eta^1(t)$ and  $\eta^2(t)$ are connected. If  $\eta(t)=\eta^1(t)\cup\eta^2(t)$ is connected , it follows from Lemma~\ref{connexdeux} that
$\eta^2(t) \cap \froufrou \eta(t)=\text{Int}(\Gamma(\eta(t))\cap\Delta_{\eta^2(t)})$ is connected.
Otherwise, $\eta^2(t) \cap \froufrou \eta(t)=\text{Int}(\eta^2(t))$, which is also connected. 
\end{proof}

We can now proceed to the proof of Theorem~\ref{th-densite}.

\begin{proof}[Proof of Theorem~\ref{th-densite}]
Let $\alpha>0$, and $\eta>0$ such that $1+\epsilon=(2+\alpha)(1/2+\eta)<1+\alpha$. Let $\alpha'>\alpha$ be such that 
$2+\alpha>1+\alpha'$ (this last condition is only necessary to ensure the good definition of the event $A_n^2$ introduced below).
For $n>0$, we define:
\begin{eqnarray*}
A_n^1  & =  & \left\{ \mathcal B_{p_1} \left( n^{2+\alpha} -4n^{1+\epsilon} \right) \subset \eta \left( n^{2+\alpha} -2n^{1+\epsilon} \right) \subset \mathcal B_{p_1} \left( n^{2+\alpha} \right) \right\}, \\
A_n^2  & = &  \Shade \left(n^{2+\alpha} -2n^{1+\epsilon}, n^{1+\alpha'} \right)^c, \\
A_n^3  & =  & \bigcap_{x \in  \mathcal B_{p_1} \left( n^{2+\alpha} \right) \backslash \mathcal B_{p_1} \left( n^{2+\alpha} -5n^{1+\epsilon} \right) } 
\left\{ \mathcal B_{p_2}^x  \left( 4(2+\alpha)n^{1+\alpha}  \right) \subset \mathcal \mathcal B_{p_2}^x  \left( 5(2+\alpha)n^{1+\alpha}  \right)\right\}, \\
A_n & = &  \bigcap_{1  \le i \le 3} A_n^i.
\end{eqnarray*}

\noindent
\underline{Step 1:}
Let us see that 
\begin{equation} \label{Borel-Cantel}
\P\left(\mathcal G^1 \cap \miniop{}{\limsup}{n\to +\infty} A_n^c\right)=0.
\end{equation}
Indeed, if $\varphi(n)=n^{2+\alpha} -2n^{1+\epsilon}$, then for $n$ large enough:
$$
n^{2+\alpha} - 4n^{1+\epsilon} \le \varphi(n) - \varphi(n)^{1/2+\eta} \le \varphi(n) + \varphi(n)^{1/2+\eta} \le n^{2+\alpha}. 
$$
By Lemma~\ref{pasdexcroissance}, there exist two strictly positive constants $A_1$ and $B_1$ such that
$$ \forall n \ge 1 \quad \P\left( \mathcal G^1 \cap (A_n^1)^c\right) \le A_1 \exp(-B_1n^{\eta (2+\alpha)}).$$
Then, as $\varphi(n)^{1/2+\eta}=o \left(n^{1+\alpha'} \right)$, by Lemma~\ref{avantage3},  there exist two strictly positive constants $A_2$ and $B_2$ such that
$$ \forall n \ge 1 \quad \P\left( \mathcal G^1 \cap  (A_n^2)^c\right) \le A_2 \exp(-B_2n^{\eta (2+\alpha)}).$$
Finally, by the large deviations result (Proposition~\ref{shapeGD}) for the $p_2$-infection,  there exist two strictly positive constants $A_3$ and $B_3$ such that
$$\forall n \ge 1 \quad  \P(  (A_n^3)^c) \le \Card{ ( \mathcal B_{p_1} \left( n^{2+\alpha} \right))} A_3 \exp(-B_3n^{1+\alpha}).$$
Collecting the previous estimates, as $\eta (2+\alpha)< 1+\alpha$, there exist two strictly positive constants $A_4$ and $B_4$ such that
\begin{eqnarray*}
\forall n \ge 1 \quad \P(\mathcal G^1 \cap A_n^c) & \le & \sum_{1  \le i \le 3} \P(\mathcal G^1 \cap (A_n^i)^c) \le A_4\exp\left(-B_4n^{\eta (2+\alpha)}\right),
\end{eqnarray*}
which proves (\ref{Borel-Cantel})
by the Borel-Cantelli lemma.

Denote by 
\begin{eqnarray*}
\Gamma_n & = & \Gamma \left( \eta \left( n^{2+\alpha} -2n^{1+\epsilon} \right) \right), \\
F_n & = & s(\Gamma_n) \cap \eta \left( n^{2+\alpha} -2n^{1+\epsilon} \right), \\
F_n^i & = & s(\Gamma_n) \cap \eta^i \left( n^{2+\alpha} -2n^{1+\epsilon} \right) 
 = \eta^i(\infty) \cap \froufrou \eta \left( n^{2+\alpha} -2n^{1+\epsilon} \right).
\end{eqnarray*}
If $n$ is large enough, $\eta \left( n^{2+\alpha} -2n^{1+\epsilon} \right) $ is connected and $\Gamma_n$ is thus a circuit. Note also that by Lemma~\ref{connexeetoile}, $F_n^2$ is $*$-connected.

\vsp
\noindent
\underline{Step 2:} 
Assume that $\mathcal G_1 \cap A_{n-1} \cap A_n$ occurs for some $n$ large enough. Let us prove that there exists $y \in \mathcal B_{p_1} \left( n^{2+\alpha}\right) \backslash  \mathcal B_{p_1} \left( n^{2+\alpha} -5n^{1+\epsilon} \right)$ such that 
\begin{eqnarray}
F_n^2& \subset & \mathcal B_{2}^y \left( 5n^{1+\alpha'}\right) \cap \mathcal B_{p_1} \left( n^{2+\alpha}\right) \backslash  \mathcal B_{p_1} \left( n^{2+\alpha} -5n^{1+\epsilon} \right). \label{6-4-step1}
\end{eqnarray} 
By $A_n^1$, $F_n \subset \mathcal B_{p_1} \left( n^{2+\alpha} \right) \backslash  \mathcal B_{p_1} \left( n^{2+\alpha} -5n^{1+\epsilon} \right)$. Assume then by contradiction that 
$$(H) \quad \exists z_1, z_2\in F_n^2 \text{ such that } \|z_1-z_2\|_2 \ge 5n^{1+\alpha'}.$$ 
As in dimension 2, 
$F_n^2$ is $*$-connected, and as $\mathcal G_1$ occurs, we can define  $\gamma_n$ as the portion of $\Gamma_n$ between $z_1$ and $z_2$ such that:
$$s(\gamma_n)\cap \eta \left( n^{2+\alpha} -2n^{1+\epsilon} \right) \subset \eta^2(\infty).$$
As the width of the annulus $\mathcal B_{p_1} \left( n^{2+\alpha} \right) \backslash  \mathcal B_{p_1} \left( n^{2+\alpha} -5n^{1+\epsilon} \right)$ is of order $5n^{1+\epsilon}=o(n^{1+\alpha'})$, there exists $\chapeau x_0 \in \mathcal S_2$ such that 
in $ \text{Cyl}_+ (\chapeau x_0, 2n^{1+\alpha'})$, the set $s(\gamma_n)\cap \eta \left( n^{2+\alpha} -2n^{1+\epsilon} \right) \subset F_n^2$ disconnects $0$ from infinity. 
By $A_n^2$, the event $\Shade(n^{2+\alpha} -2n^{1+\epsilon}, 2n^{1+\alpha'})$ can not occur, thus in $ \text{Cyl}_+ (\chapeau x_0, 2n^{1+\alpha'})$, the set $F_n^2$ does not disconnect $\eta^1(n^{2+\alpha} -2n^{1+\epsilon})$ from infinity. Let $z \in \eta^1(n^{2+\alpha} -2n^{1+\epsilon})$ be in the same connected component as infinity in $\text{Cyl}_+ (\chapeau x_0, 2n^{1+\alpha'})$ deprived of $F_n^2$. Note that -- this is a key point -- the infection path from $s_1$ to $z$ has to enter
$\text{Cyl}_+ (\chapeau x_0, 2n^{1+\alpha'}) \cap \mathcal B_{p_1} \left( n^{2+\alpha} \right) \backslash  \mathcal B_{p_1} \left( n^{2+\alpha} -5n^{1+\epsilon} \right)$ by crossing the border of the cylinder.

If $n$ is large enough, then $(n-1)^{2+\alpha}<n^{2+\alpha}-5n^{1+\epsilon}$. By $A_{n-1}^1$, the set $\eta \left( (n-1)^{2+\alpha} -2(n-1)^{1+\epsilon} \right)$ is contained in $\mathcal B_{p_1} \left( (n-1)^{2+\alpha} \right)$, and thus, an infection path from 
$s_1$ to $z$ has to visit some vertex $s \in \eta^1(\infty)$ satisfying $s \in \mathcal B_{p_1} \left( n^{2+\alpha} -5n^{1+\epsilon} \right) \backslash \mathcal B_{p_1} \left( (n-1)^{2+\alpha} \right)$, and thus such that 
$$s \not\in B^{s_1}_{p_1}((n-1)^{2+\alpha} -2(n-1)^{1+\epsilon}).$$
This implies 
$$z \in B_{p_1}^{s}(n^{2+\alpha} -2n^{1+\epsilon} -(n-1)^{2+\alpha}+2(n-1)^{1+\epsilon}) \subset B_{p_1}^{s}(2(2+\alpha)n^{1+\alpha}).$$
Thus, $z$ has to be at a random distance for the $p_2$-infection less than $2(2+\alpha)n^{1+\alpha}$ of a point in the border of the cylinder, and, by $A_n^3$ at a distance $\norme{.}{p_2}$ less than $5(2+\alpha)n^{1+\alpha}$ of a point in the border of the cylinder; and this is also the case for any point $z \in  \eta^1(n^{2+\alpha} -2n^{1+\epsilon})$ in the same connected component as infinity in $\text{Cyl}_+ (\chapeau x_0, 2n^{1+\alpha'})$ deprived of $F_n^2$. 
As $1+\alpha'>1+\alpha$, this implies, for $n$ large enough, that in $ \text{Cyl}_+ (\chapeau x_0, n^{1+\alpha'})$, the set $F_n^2$ disconnects $\eta^1(n^{2+\alpha} -2n^{1+\epsilon})$ from infinity, which contradicts $A_n^2$, and thus $(H)$. This completes the proof of (\ref{6-4-step1}).

\vsp
\noindent
\underline{Step 3:}
Assume now that $n$ is large enough and that $\mathcal G_1 \cap A_{n-1} \cap A_n \cap A_{n+1} \cap A_{n+2}$ occurs. Let $y \in \mathcal B_{p_1} \left( n^{2+\alpha} \right) \backslash  \mathcal B_{p_1} \left( n^{2+\alpha} -5n^{1+\epsilon} \right)$ be such that (\ref{6-4-step1}) is satisfied. Let us prove that
\begin{equation}
\label{6-4-step2}
\eta^2(\infty) \cap \left(  \mathcal B_{p_1} \left( (n+1)^{2+\alpha} \right) \backslash \mathcal B_{p_1} \left( n^{2+\alpha}  \right) \right) \subset \mathcal B_2^{y} (6n^{1+\alpha'}).
\end{equation}
Consider
$z \in \eta^2(\infty) \cap \left(  \mathcal B_{p_1} \left( (n+1)^{2+\alpha} \right) \backslash \mathcal B_{p_1} \left( n^{2+\alpha}  \right) \right).$
If $n$ is large enough, $(n+1)^{2+\alpha} \le (n+2)^{2+\alpha} -4(n+2)^{1+\epsilon}$, and thus, by $A_{n+2}^1$, the infection time of $z$ is less than $(n+2)^{2+\alpha} -2 (n+2)^{1+ \epsilon}.$
By $A_{n}^1$, as $z \not\in \mathcal B_{p_1} \left( n^{2+\alpha} \right)$, its infection time is also strictly larger than $n^{2+\alpha} -2 n^{1+ \epsilon}.$ Thus the infection path from $s_2$ to $z$ has to visit some point $s \in F_n^2$.
By (\ref{6-4-step1}), the point $ s$ is not in $\mathcal B_{p_1} \left( n^{2+\alpha} -5n^{1+\epsilon} \right)$; thus, as $(n-1)^{2+\alpha} <n^{2+\alpha} -5n^{1+\epsilon}$ for $n$ large enough, $A_{n-1}^1$ ensures that the infection time for $s$ is larger than $ (n-1)^{2+\alpha} -2(n-1)^{1+\epsilon}$, which leads, for $n$ large enough, to
$$ z \in B_{p_2}^s((n+2)^{2+\alpha} +2 (n+2)^{1+ \epsilon}-( (n-1)^{2+\alpha} -2(n-1)^{1+\epsilon}))
 \subset B_{p_2}^s(4(2+\alpha)n^{1+\alpha}).$$
As $s \in \mathcal B_{p_1} \left( n^{2+\alpha} \right) \backslash \mathcal B_{p_1} \left( n^{2+\alpha} -5n^{1+\epsilon} \right)$, by $A_n^3$,
we have $z \in \mathcal B_{p_2}^s(5(2+\alpha)n^{1+\alpha})),$ and thus
$$z \in \bigcup_{s \in \mathcal B_2^{y} (5n^{1+\alpha'}) } \mathcal B_{p_2}^s(5(2+\alpha)n^{1+\alpha})),$$
which is included in $\mathcal B_2^{y} (6n^{1+\alpha'})$ for $n$ large enough. This proves (\ref{6-4-step2}).

\vsp
\noindent
\underline{Step 4:} Let us prove that there exists a constant $C$ such that 
\begin{equation}
\label{6-4-step3}
\forall n \ge 1 \quad | \eta^2(\infty) \cap \mathcal B_{p_1}(n^{2+\alpha})|
\le C n^{3+2 \alpha'}.
\end{equation}
By (\ref{Borel-Cantel}), on $\mathcal G^1$, there almost surely exists $m$ such that $A_n$ occurs for every $n \ge m$. Thus, almost surely, for every $n > m$, by (\ref{6-4-step1}) and (\ref{6-4-step2}) -- and by increasing $m$ if necessary --, there exists $y_n \in \mathcal B_{p_1} \left( n^{2+\alpha} \right) \backslash  \mathcal B_{p_1} \left( n^{2+\alpha} -5n^{1+\epsilon} \right)$ such that
$$
\eta^2(\infty) \cap \left(  \mathcal B_{p_1} \left( (n+1)^{2+\alpha} \right) \backslash \mathcal B_{p_1} \left( n^{2+\alpha}  \right) \right) \subset \mathcal B_2^{y_n} (6 n^{1+\alpha'}).
$$
Thus, there exist constants $C_i$ such that for $n \ge m$,
\begin{eqnarray*}
| \eta^2(\infty) \cap \mathcal B_{p_1}(n^{2+\alpha})| & \le & C_1 \Card{\mathcal B_{p_1}(m^{2+\alpha})} +\sum_{k=m}^{n-1} C_2 \Card{\mathcal B_2 (6k^{1+\alpha'})} \\
& \le & C_3 m^{4+2 \alpha} + C_4 \sum_{k=m}^{n-1} k^{2+2 \alpha'} \le  C_5 n^{3+2 \alpha'}.
\end{eqnarray*}

\vsp
\noindent
\underline{Step 5:} 
Consider $\beta \in (0,1/2)$. We can choose $\alpha'>\alpha>0$ such that $\frac{3+2\alpha'}{2+\alpha}=3/2+\beta$ and such that $1+\alpha'<2+\alpha$.
For $t$ large enough, choose $n$ such that
$n^{2+\alpha} \le t <(n+1)^{2+\alpha}.$
Thus, by the previous step,
$$
| \eta^2(\infty) \cap \mathcal B_{p_1}(t)|  \le  | \eta^2(\infty) \cap \mathcal B_{p_1}( (n+1)^{2+\alpha}) | \le C (n+1)^{3+2 \alpha'}\sim  Ct^{\frac{3+2\alpha'}{2+\alpha}}=Ct^{3/2+\beta}.
$$
The norm equivalence implies the analogous result for any norm on $\R^2$, which proves the first point of Theorem~\ref{th-densite}.

Let us now  prove the second point:
for every $\beta>0$, 
$$\miniop{}{\limsup}{t\to +\infty} \frac{ \text{Diam}\left(\left(\eta^2(\infty)+[-1/2,1/2]^2 \right) \cap \mathcal S_{p_1}(t)\right)}{t^{1/2+\beta}}=0.$$
It is clearly sufficient to consider $\beta \in (0,1/2)$.
We can choose $\alpha'>\alpha>0$ such that $\frac{1+\alpha'}{2+\alpha}<1/2+\beta$.
For $t$ large enough, choose $n$ such that
$n^{2+\alpha} < t \le (n+1)^{2+\alpha}.$
Then, by  (\ref{6-4-step2}), there exists $y \in \R^2$ such that
$$\left(\eta^2(\infty)+[-1/2,1/2]^2 \right) \cap \mathcal S_{p_1}(t)\cap\Z^2 \subset \mathcal B_{p_1} ( (n+1)^{2+\alpha})  \backslash \mathcal B_{p_1} ( n^{2+\alpha} ) \subset \mathcal B_2^{y} (6n^{1+\alpha'}).$$
As $n \sim t^{\frac1{{2+\alpha}}}$, this ends the proof of the second point.

Turning to the proof of the last assertion, consider the following
alternative:
\begin{itemize}
\item If the weak species (type~1) does not unboundedly grow, its natural density is zero, while the density of the strong is one.
\item If the weak species grows unboundedly, the first point of the present
theorem ensures that the strong species has null density, and therefore that
the weak have full density.
\end{itemize}
\end{proof}

\section{Non-coexistence except perhaps for a denumerable set}

In this section, we prove Theorem~\ref{HPgeneral}: Remember that H\"aggstr\"om and Pemantle proved non-coexistence for two epidemics progressing according exponential laws with parameter $1$ and $\lambda$ "except perhaps for a denumerable set" for $\lambda$, and Theorem~\ref{HPgeneral} extends this result to families of laws depending on a continuous parameter.

The first step consists in coupling all possible competition models on the same probability space, respecting the stochastic order of the laws. This will give natural inclusions between sets of infected points for competition with distinct parameters, as stated in Lemma~\ref{yal}. Assume that $p<r$, that the slow infection (the strong one) uses the law with parameter $p$ (respectively,  $r$) and that both infections manage to grow unboundedly. Let also choose $q\in (p,r)$.
Then, we can expect than strengthening the slow infection by increasing its parameter from $p$ to $q$ makes it strong enough to win and surround the fast one. Proving this is the aim of Lemma~\ref{yaol}, and the proof is based on Theorem~\ref{pasdexcroissance}. Finally, to prove Theorem~\ref{HPgeneral}, we show that for a fixed $q \in I$,  the set of $p<q$ such that $\P(\mathcal{G}^1_{p, q}\cap \mathcal{G}^2_{p, q})>0$ is a subset of the discontinuity set of an increasing function, and is thus at most denumerable. 

\subsection*{Coupling} 
We first couple all passage times for varying parameters
thanks to the generalized inverse of the repartition function, 
which allows to build all the competition models 
 on the same probability space. 
This generalizes the construction presented in the introduction with only
two parameters.

On $\Omega=[0,1]^{\Ed}$, consider the probability measure $\P=\mathcal U[0,1]^{\otimes\Ed}$, where $\mathcal U[0,1]$ denotes the uniform law on the set $[0,1]$. 
For each $\omega \in\Omega$ and $p\in I$, define
$$t^p_e=\inf\{x\in\R: \;\nu_p((-\infty,x])\ge \omega_e\}.$$
Under $\P$, the variables $(t^p_e)_{e\in\Ed}$  are independent identically distributed  with common law $\nu_{p}$.
Moreover 
$$\forall e\in\Ed \quad \forall (p,q)\in I^2 \quad p\le q \Longrightarrow t^p_e \ge t^q_e.$$

We build now, for a given $(p_1,p_2)\in I^2$, the competition process in a realization $\omega \in \Omega$. 
Let  $E=([0,+\infty]\times [0,+\infty] )^{\Zd}$.
We recursively define a $E$-valued sequence $(X_n^{p_1,p_2})_{n\ge 0}$ and a non-negative
sequence  $(T_n^{p_1,p_2})_{n\ge 0}$. The sequence $(T_n^{p_1,p_2})_{n\ge 0}$ contains the successive times of infections, while a point $\epsilon=(\epsilon^{1}(z), \epsilon^2(z))_{z \in \Zd} \in E$ codes, for each site $z$, its times of infection $\epsilon^1(z)$ ( $\epsilon^2(z)$)  by the first (respectively, the second) infection. We start the process with two distinct sources $s_1$ and $s_2$ in $\Zd$, and set $T_0^{p_1,p_2}=0$ and 
$$
X_0^{p_1,p_2}=(X_0^{1,p_1,p_2}(z), X_0^{2,p_1,p_2}(z))_{z \in \Zd} \text{ with } 
\left\{
\begin{array}{l}
X_0^{i,p_1,p_2}(z)=0 \text{ if } z=s_i, \\
X_0^{i,p_1,p_2}(z)=+\infty \text{ otherwise.}
\end{array}
\right.
$$
This means that at time $0$, no point of $\Zd$ has been infected yet but the two initial sources $s_1$ and $s_2$. Then, for $n \ge 0$, define the next time of infection:
$$T_{n+1}^{p_1,p_2}=\inf\{ X_n^{i,p_1,p_2}(y)+t_{\{y,z\}}^{p_i}: \quad \{y,z\}\in\Ed, \; i\in\{1,2\}, \;  X_n^{3-i,p_1,p_2}(z)=+\infty\}.$$
Note that the infimum in the definition of $T_{n+1}$ is always taken on a finite set. Moreover,
Assumption ($5$) ensures that if this infimum is reached by several triplets $(i,y,z)$, all these triplets have the same first coordinate, which means that a point can be infected by the same species from  distinct neighbors at the same time, but not by the two species simultaneously. For such a triplet, the next infection is of type $i$ from (one of the) $y$ to $z$. The set of infected points of type $3-i$ has not changed:
$$\forall x \in \Zd \quad  X_{n+1}^{3-i,p_1,p_2}(x)=X_n^{3-i,p_1,p_2}(x),$$
while the point $z$ has been infected by species $i$ at time $X_n^{i,p_1,p_2}(y)+t_{\{y,z\}}^{p_i}$:
$$\forall x \in \Zd\backslash \{z\} \quad X_{n+1}^{i,p_1,p_2}(x)=X_n^{i,p_1,p_2}(x) \text{ and } X_{n+1}^{i,p_1,p_2}(z)=X_n^{i,p_1,p_2}(y)+t_{\{y,z\}}^{p_i}.$$
Note that $X_n^{i,p_1,p_2}(y)$ and $X_n^{3-i,p_1,p_2}(y)$ can not be simultaneously finite, which corresponds to the fact that each site is infected by at most one type of infection. Moreover, once $\min(X_n^{1,p_1,p_2}(x), X_n^{2,p_1,p_2}(x))$ is finite, its value -- the time of infection of $x$ -- does not change any more.

As $\nu_{p_1} \gstoch \nu_{p_2}$, species $1$ is slower than species $2$. We also define the sets $\eta^{p_1,p_2}(t),\eta^{1,p_1,p_2}(t),\eta^{2,p_1,p_2}(t)$ that are respectively infected points, infected points of type $1$, infected points of type $2$ at time $t$, by
$$\forall i\in\{1,2\}\quad \forall t\in [T_n^{p_1,p_2},T_{n+1}^{p_1,p_2})\quad  \eta^{i,p_1,p_2}(t)=\{z\in\Zd; X_n^{i,p_1,p_2}(z)<+\infty\}$$
and $\eta^{p_1,p_2}(t)=\eta^{1,p_1,p_2}(t)\cup \eta^{2,p_1,p_2}(t)$. 
We also define
\begin{eqnarray*}
\eta^{i,p_1,p_2}(\infty) & = & \miniop{}{\cup}{t\ge 0}\eta^{i,p_1,p_2}(t), \\
\mathcal{G}^{i,p_1,p_2} & = & \left\{ |\eta^{i,p_1,p_2}(\infty)|=+\infty\right\} \text{ for } i=1,2.
\end{eqnarray*}
The set $\mathcal{G}^{i,p_1,p_2}$ corresponds to the survival of type $i$.

\begin{lemme} 
\label{yal}
Let $t\ge 0$.
\begin{itemize}
\item $\eta^{1,p,q}(t)$ is non-decreasing in $p$ and non-increasing in $q$,
\item $\eta^{2,p,q}(t)$ is non-decreasing in $q$ and non-increasing in $p$.
\end{itemize}
\end{lemme}
\begin{proof} The proof is just a "several parameters" version of the proof of Lemma~\ref{coupling}. 
We only prove the monotonicity of the two sets with respect to $q$.
Let  then $q<r$. We now proceed by induction to prove that for every $n \in \N$ 
$$ (H_n) \quad \forall x \in \Zd \quad X_n^{1,p,q}(x)\le X_n^{1,p,r}(x) \text{ and }X_n^{2,p,q}(x)\ge X_n^{2,p,r}(x).$$
Clearly, $(H_0)$ is true. Assume that $(H_n)$ holds. 

1. Let us first prove that $X_{n+1}^{1,p,q}(x)\le X_{n+1}^{1,p,r}(x)$
for each $x$ in $\Zd$. We have the following alternative:
\begin{itemize}
\item If $X_{n+1}^{1,p,r}(x)=+\infty$, it is obvious that $X_{n+1}^{1,p,q}(x)\le X_{n+1}^{1,p,r}(x).$
\item If $X_{n}^{1,p,r}(x)<+\infty$ then, by $(H_n)$, $X_{n}^{1,p,q}(x)$ is also finite, 
and thus their values do not change when we go from $n$ to $n+1$, and the inequality is preserved.
\item If $X_{n+1}^{1,p,r}(x)<+\infty$ and $X_{n}^{1,p,r}(x)=+\infty$, there exists $y \sim x$ which has infected $x$ at the $(n+1)$-th step of the construction with type $1$ species. In other words:
$$X_{n+1}^{1,p,r}(x)=X_{n}^{1,p,r}(y)+t^p_{\{y,x\}} \text{ and } \forall z \sim x, \; X_{n}^{1,p,r}(y)+t^p_{\{y,x\}}<  X_{n}^{2,p,r}(z)+t^r_{\{z,x\}}.$$
By $(H_n)$, since $t^r\le t^q$, we have for each $z \sim x$:
$$X_{n}^{1,p,q}(y)+t^p_{\{y,x\}} \le X_{n}^{1,p,r}(y)+t^p_{\{y,x\}} <  X_{n}^{2,p,r}(z)+t^r_{\{z,x\}} \le X_{n}^{2,p,q}(z)+t^q_{\{z,x\}}.$$
This says that, in the $(p,q)$ competition, $x$ is infected by the species $1$ and that $X_{n+1}^{1,p,r}(x)=X_{n}^{1,p,r}(y)+t^p_{\{y,x\}} \ge X_{n+1}^{1,p,q}(x)$.
\end{itemize}

2. Let us now prove that $X_{n+1}^{2,p,q}(x)\ge X_{n+1}^{2,p,r}(x)$
for each $x \in \Zd$. We have the following alternative:
\begin{itemize}
\item If $X_{n+1}^{2,p,q}(x)=+\infty$, it is obvious that $X_{n+1}^{2,p,q}(x)\ge X_{n+1}^{2,p,r}(x).$
\item If $X_{n}^{2,p,q}(x)<+\infty$ then, by $(H_n)$, $X_{n}^{2,p,r}(x)$ is also finite, 
and thus their values do not change when we go from $n$ to $n+1$, and the inequality is preserved.
\item If $X_{n+1}^{2,p,q}(x)<+\infty$ and $X_{n}^{2,p,q}(x)=+\infty$, there exists $y \sim x$ which has infected $x$ at the $(n+1)$-th step of the construction with type $2$ species. In other words:
$$X_{n+1}^{2,p,q}(x)=X_{n}^{2,p,q}(y)+t^q_{\{y,x\}} \text{ and } \forall z \sim x, \; X_{n}^{2,p,q}(y)+t^q_{\{y,x\}}<  X_{n}^{1,p,q}(z)+t^p_{\{z,x\}}.$$
By $(H_n)$, since $t^r\le t^q$, we have for each $z \sim x$:
$$X_{n}^{2,p,r}(y)+t^r_{\{y,x\}} \le X_{n}^{2,p,q}(y)+t^q_{\{y,x\}} < X_{n}^{1,p,q}(z)+t^p_{\{z,x\}} \le X_{n}^{1,p,r}(z)+t^p_{\{z,x\}}.$$
This says that, in the $(p,r)$ competition, $x$ is infected by the species $2$ and that $X_{n+1}^{2,p,r}(x) \le X_{n}^{2,p,r}(y)+t^r_{\{y,x\}} =X_{n+1}^{2,p,q}(x)$.
\end{itemize} 
To conclude, note that
$\eta^{i,p,q}(t)  =  \{ z \in \Zd: \; \exists n \in \N, \; X^{i,p,q}(z) \le t \}.$
It is then obvious that 
$\eta^{1,p,q}(t)$ is non-increasing in $q$, and 
$\eta^{2,p,q}(t)$ is non-decreasing in $q$.
\end{proof}

\begin{lemme}
\label{shape}
For $A\subset\Zd$, we define
$$|A|_p=\sup\{\norme{x}{p}: \; x\in A\} \quad
\text{ and }
\quad |A|_{*,p}=\inf\{\norme{x}{p}: \; x\in \Z^d \backslash A\}.$$
Then, $\P$-almost surely,
$$\frac{|B^a_{p}(t)|_{p}}t\to 1 \quad 
\text{ and }
\quad \frac{|B^a_{p}(t)|_{*,p}}t\to 1.$$
\end{lemme}
\begin{proof}
These are direct consequences of the large deviation result, Proposition~\ref{shapeGD}.
\end{proof}

\begin{lemme}
\label{yaol}
If $p<q<r$, then
$\P(\mathcal{G}^{1,p,r}\cap \mathcal{G}^{2,q,r})=0$.
\end{lemme}

\begin{proof}
By coupling, $\mathcal{G}^{2,p,r}\supset \mathcal{G}^{2,q,r}$, thus we have
$\mathcal{G}^{1,p,r}\cap \mathcal{G}^{2,q,r}=(\mathcal{G}^{1,p,r}\cap \mathcal{G}^{2,p,r})\cap \mathcal{G}^{2,q,r}$.
So we can assume that $\mathcal{G}^{1,p,r}\cap \mathcal{G}^{2,p,r}$ occurs and prove that $\mathcal{G}^{2,q,r}$ can not happen.
By Theorem~\ref{pasdexcroissance}, we have 
$$
\miniop{}{\limsup}{t\to +\infty}\frac{|\eta^{2,p,r}(t)|_{p}}t\le 1, 
\text{ which implies, by Lemma~\ref{yal}, }
\miniop{}{\limsup}{t\to +\infty}\frac{|\eta^{2,q,r}(t)|_{p}}t\le 1.
$$
Now, by Proposition~\ref{strictecomp}, we have
$$\miniop{}{\limsup}{t\to +\infty}\frac{|\eta^{2,q,r}(t)|_{q}}t\le C_{p,q}<1.$$
Using the coupling Lemma~\ref{coupling} and Lemma~\ref{shape} together, we get
$$\miniop{}{\liminf}{t\to +\infty}\frac{|\eta^{1,q,r}(t)\cup \eta^{2,q,r}(t)|_{*,q}}t\ge 1.$$
Now, let $t$ be large enough to ensure that 
$$\frac{|\eta^{2,q,r}(t)|_{q}}t\le\frac{2C_{p,q}+1}3=\alpha \;
\text{ and } \;
\frac{|\eta^{1,q,r}(t)\cup \eta^{2,q,r}(t)|_{*,q}}t\ge \frac{C_{p,q}+2}3=\beta.$$
Then every point $x$ such that $\alpha t <\|x\|_{q}< \beta t$
belongs to $\eta^{1,q,r}(t)\backslash\eta^{2,q,r}(t)$, which prevents the occurrence of the event
$\mathcal{G}^{2,q,r}$.
\end{proof}

\begin{proof}[Proof of Theorem~\ref{HPgeneral}]
Let $q \in I$ be fixed and consider the maps
$s^i:p\mapsto \P(\mathcal{G}^{i,p, q})$.
By Lemma~\ref{yal}, $s^1$ is  non-decreasing, whereas $s^2$ is  non-increasing.
Suppose now that $p<q$.
We prove that $\P(\mathcal{G}^{1,p, q}\cap \mathcal{G}^{2,p, q})=0$
if $s^1$ 
is left-continuous 
at $p$. By Lemma~\ref{yaol}, 
\begin{eqnarray*}
\P(\mathcal{G}^{1,p, q}\cap \mathcal{G}^{2,p, q}) & = & 
\P(\mathcal{G}^{1,p-1/n, q}\cap\mathcal{G}^{1,p, q}\cap \mathcal{G}^{2,p, q})+
\P((\mathcal{G}^{1,p-1/n, q})^c\cap\mathcal{G}^{1,p, q}\cap \mathcal{G}^{2,p, q})\\
& = & \P((\mathcal{G}^{1,p-1/n, q})^c\cap\mathcal{G}^{1,p, q}\cap \mathcal{G}^{2,p, q})\\
& \le & \P((\mathcal{G}^{1,p-1/n, q})^c\cap\mathcal{G}^{1,p, q})\\
& \le & \P(\mathcal{G}^{1,p, q})-\P(\mathcal{G}^{1,p-1/n, q})\\
& \le & s^1(p)- s^1(p-1/n).
\end{eqnarray*}
Thus, the set of $p<q$ such that $\P(\mathcal{G}^1_{p, q}\cap \mathcal{G}^2_{p, q})>0$ is a subset of the discontinuity set of the non-decreasing function $s^1$. Therefore, it is at most denumerable.
Note that we could prove that $\P(\mathcal{G}^{1,p, q}\cap \mathcal{G}^{2,p, q})=0$ if $s^2$ is  right-continuous at $p$ in the very same way. 
\end{proof}
\section{Concluding remarks}
Since we are coming to the end of our study, it is worth questioning
 the relevance of the notion of strong coexistence and its relationship
with the Häggström-Pemantle conjecture.

At first, let us say a word about the notion of strong coexistence.
It is obviously stronger, but how far away is it from  the notion of coexistence? In the case where the two species have the same passage times law,
a partial answer is given by a recent work by Gou\'er\'e~\cite{gouere}: 
some of these results imply that -- under classical assumptions implying coexistence -- one can find initial configurations that give rise to strong coexistence with positive probability. The restriction on the initial conditions can however be dropped by a modification argument as in Garet-Marchand~\cite{GM-coex}.
Actually, in the case where the two species have the same passage times law,
simulations let think that when coexistence occurs, each species grants itself
a cone, and thus strong coexistence occurs.

The above remarks seem to show the relevance of the notion of strong coexistence, but we must now wonder how far we are from the Häggström-Pemantle conjecture.
First, it could be interesting to note that Theorem~\ref{HPgeneral} allows a reformulation of this conjecture: it is sufficient to prove that the map
$p\mapsto \P_{p,q}(\mathcal G^1)$ has a unique point of discontinuity.
Note also that this reformulation is quite close to some recent result by
Deijfen and Häggström concerning competition with exponential speeds on general graphs (Theorem 4.1 in~\cite{DH-nonmon}). However, it is not evident that this remark can be exploited to prove the conjecture.
In this perspective, it it more natural to look at Theorem~\ref{avantage3} and Corollary~\ref{avantage4}: it seems highly unlikely that the strong could 
survive being constrained to occupy only a negligible portion of the aerial 
surface, but we did not success to prove it at this time.


\def\refname{References}
\bibliographystyle{plain}
\bibliography{rikiki}


\end{document}